# A local discontinuous Galerkin gradient discretization method for linear and quasilinear elliptic equations

Assyr Abdulle[*], Giacomo Rosilho de Souza[†]

July 27, 2018

**Abstract**

A local weighted discontinuous Galerkin gradient discretization method for solving elliptic equations is introduced. The local scheme is based on a coarse grid and successively improves the solution solving a sequence of local elliptic problems in high gradient regions. Using the gradient discretization framework we prove convergence of the scheme for linear and quasilinear equations under minimal regularity assumptions. The error due to artificial boundary conditions is also analyzed, shown to be of higher order and shown to depend only locally on the regularity of the solution. Numerical experiments illustrate our theoretical findings and the local method's accuracy is compared against the non local approach.

## 1 Introduction

Partial differential equations with high contrast are notoriously difficult to solve. In order to capture strong variations of the exact solution in the numerical approximations of the PDE, non uniform grids are usually required. The construction of such grids is often based on an iterative process, where a solution is computed and an a posteriori error estimator is used to indicate the regions where the mesh has to be refined, see [2, 3, 21, 23]. In such approach, the solution is computed on the whole domain at each step, even if the mesh has changed only in a small portion of the domain.

In this paper we propose an algorithm for elliptic PDE, based on a decomposition of the computational domain in local subdomains adapted to the variation of the solution. In each subdomain, only local problems need to be solved and no iterations are needed between subdomains (e.g. as in domain decomposition method), as we define artificial boundary conditions and do compute the solution only once in each local domains. We concentrate here on the a priori error analysis of our scheme, while we postpone the a posteriori error analysis to a companion paper [1]. The local scheme proposed in this paper is more efficient than the classical schemes for elliptic PDEs with strong variations for several reasons.

- For linear problems, when using an iterative solver such as the conjugate gradient (CG) method we have smaller problems to compute on the finer meshes, while the non-local classical schemes need the solution of global linear systems with a large number of degrees of freedom (DOF) (recall that the CG method has a convergence rate that is super-linear with respect to the DOF of the system).

---

[*]École Polytechnique Fédérale de Lausanne, Switzerland (assyr.abdulle@epfl.ch)

[†]École Polytechnique Fédérale de Lausanne, Switzerland(giacomo.rosilhodesouza@epfl.ch).

- When solving a linear system arising from PDEs with CG methods, preconditioners are ususally needed, a usual choice for CG being the incomplete Cholesky (IC) factorization. For non-local schemes, the high contrast of the PDE leads to systems with high condition number (due to mesh and data variations). For the local scheme, as each subdomain involves smaller variations of the solution the condition number is smaller, leading to faster convergence of the iterative method.

- Finally for nonlinear problems, in addition to the computational saving described previously, only a nonlinear problem on a coarse global mesh needs to be solved for the local scheme, while the subsequent local problems are linear. The computational saving is therefore significant for such problems.

The idea of solving local elliptic problems to improve the numerical solution's accuracy is not new in the literature. The Local Defect Correction method (LDC) first presented in [14] is an iterative process that at each step solves a global problem on a coarse mesh and a local problem on a fine mesh. The solution of the global problem provides artificial boundary conditions to the local problem. The solution of the latter is then introduced into the coarse system to estimate its residual. The coarse system is solved again but adding the residual to its right hand side, leading to a more accurate coarse solution and hence better artificial boundary conditions for the next local problem. Two similar methods are the Fast Adaptive Composite grid algorithm [18] and the Multi-Level Adaptive Technique [5]. In [12] it is shown that under reasonable assumptions the three methods lead to the same solution. In their original form the schemes were defined for finite difference methods but finite volumes or finite element versions exist, see [19, 24]. Only recently has the LDC scheme been coupled with an a posteriori error estimator in order to automatically select the local domains [4].

The local method that we propose in this paper relies on the discontinuous Galerkin discretization, more precisely on the Symmetric Weighted Interior Penalty Galerkin (SWIPG) scheme [10, 7]. We consider the elliptic model problem

$$-\nabla \cdot (A\nabla u) = f \quad \text{in } \Omega, \tag{1.1a}$$
$$u = 0 \quad \text{in } \partial\Omega, \tag{1.1b}$$

where $\Omega \subset \mathbb{R}^d$ for $d = 1, 2, 3$ is an open bounded polytopal connected set, $u \in H_0^1(\Omega)$ and $f \in H^{-1}(\Omega)$. The matrix $A$ is symmetric, positive definite and can possibly depend on $u$, since we consider both a linear and a quasilinear case.

The scheme that we propose is different from the aforementioned methods in the sense that it computes only one global solution on the full domain while all the subsequent computations are local. Additionally, the a priori error analysis is performed under minimal regularity assumptions, that is, assuming $u \in H_0^1(\Omega)$ and $f \in H^{-1}(\Omega)$. This is achieved by recasting the SWIPG scheme into the Gradient Discretization (GD) framework [9, 11]. The GD method is a framework suitable for studying the a priori convergence of various types of diffusion problems: linear and non linear, steady state or transient. For our scheme, the GD framework is convenient to decompose the sources of errors in the local problems. Furthermore, applying the pointwise estimates from [6], we can prove (in some particular cases) that the errors coming from the artificial boundary conditions are of higher order and depend only locally on the regularity of the solution. Finally, we stress out that the GD framework is only used for the analysis, indeed another advantage of the scheme is that it fits very easily in existing codes that use the popular discontinuous Galerkin scheme without needing additional data structures nor additional memory requirements.

The paper is organized as follows. In Section 2 we present the Gradient Discretization framework and the Symmetric Weighted Discontinuous Galerkin Gradient Discretization (SWDGGD), which is equivalent to the SWIPG scheme. At the end of the section we introduce a local version of the SWDGGD. In Section 3 we present the local scheme and establish an a priori error analysis for linear equations. In Section 4 we introduce the scheme and the a priori error analysis for quasilinear equations. Finally, Section 5 provides numerical results and comparison with the classical scheme. The equivalence between the SWDGGD and SWIPG methods is postponed to Appendix A.

## 2 Notation and preliminary results

Our local scheme is based on the traditional SWIPG scheme but the analysis is done in the GD framework, this allows for minimal assumptions and further generalizations as quasilinear problems. Whence we introduce in Section 2.1 the notation for the GD setting and in Section 2.2 we define a particular GD scheme which is equivalent to SWIPG, their equivalence is shown in Appendix A. The method presented in Section 2.2 is a slight modification of the one proposed in [11], the main difference is that the latter is equivalent to the Symmetric Interior Penalty Galerkin (SIPG) method. We opted for the SWIPG scheme instead of SIPG since it is known to have improved stability in problems with high diffusivity contrasts [10] and also to be suitable for a locally vanishing diffusion [8]. In Sections 2.1 and 2.2 we mainly follow [9] and [11]. In what follows we make the following assumptions on the data for the linear case

**Assumption 2.1.**

- $\Omega \subset \mathbb{R}^d$ is an open bounded polytopal domain,
- $A : \Omega \to \mathbb{R}^{d \times d}$ is such that $A(\boldsymbol{x})$ is a symmetric matrix measurable with respect to $\boldsymbol{x}$ and there exists $\underline{\lambda}, \overline{\lambda} > 0$ such that it has eigenvalues in $[\underline{\lambda}, \overline{\lambda}]$,
- the forcing term is $f \in H^{-1}(\Omega)$.

For the quasilinear case, we will assume

**Assumption 2.2.**

- $\Omega \subset \mathbb{R}^d$ is an open bounded polytopal domain,
- $A(\boldsymbol{x}, s) = (a_{ij}(\boldsymbol{x}, s))_{i,j=1}^d$ is such that $a_{ij} : \overline{\Omega} \times \mathbb{R} \to \mathbb{R}$ is continuous in $\boldsymbol{x}$ and Lipschitz continuous in $s$. Furthermore $A(\boldsymbol{x}, s)$ is a symmetric matrix with eigenvalues in $[\underline{\lambda}, \overline{\lambda}]$,
- the forcing term is $f \in H^{-1}(\Omega)$.

For simplicity, the dependence of $A$ on $\boldsymbol{x}$ is left out in our notation. Under Assumption 2.1 the weak solution of Eq. (1.1) is $u \in H_0^1(\Omega)$ such that

$$\int_\Omega A \nabla u \cdot \nabla v \, d\boldsymbol{x} = \langle f, v \rangle \qquad \text{for all } v \in H_0^1(\Omega), \tag{2.1}$$

where $\langle \cdot, \cdot \rangle$ denotes the pairing between $H^{-1}(\Omega)$ and $H_0^1(\Omega)$. Under Assumption 2.2 we have a similar weak solution obtained by replacing $A$ by $A(u)$ in Eq. (2.1).

## 2.1 The Gradient Discretizazion method

We start by defining the GD method for homogeneous Dirichlet boundary conditions as introduced in [11] along with some of its properties.

**Definition 2.3.** *A gradient discretization method $\mathcal{D}$ for homogeneous Dirichlet boundary conditions is defined by $\mathcal{D} = (X_\mathcal{D}, \Pi_\mathcal{D}, \nabla_\mathcal{D})$, where*

1. *the set $X_\mathcal{D}$ is a finite dimensional real vector space,*

2. *the reconstruction function $\Pi_\mathcal{D} : X_\mathcal{D} \to L^2(\Omega)$ is a linear mapping that reconstructs, from an element in $X_\mathcal{D}$, a function over $\Omega$,*

3. *the gradient reconstruction $\nabla_\mathcal{D} : X_\mathcal{D} \to L^2(\Omega)^d$ is a linear mapping which reconstructs, from an element of $X_\mathcal{D}$, a gradient over $\Omega$. This gradient reconstruction must be chosen such that $\|\nabla_\mathcal{D} \cdot \|_{L^2(\Omega)^d}$ is a norm on $X_\mathcal{D}$.*

**Example 2.4.** *Among others, the conforming $\mathbb{P}_1$ finite element Galerkin method can be written as a GD method. Given a partition of $\Omega$ into simplices, let $V_h \subset H^1_0(\Omega)$ be the set of piecewise linear and continuous functions on this partition. Let $\{e_i\}_{i \in I}$ be a basis of $V_h$, we define $X_\mathcal{D} = \{\phi = (\zeta_i)_{i \in I} : \zeta_i \in \mathbb{R} \text{ for all } i \in I\}$, $\Pi_\mathcal{D} \phi = \sum_{i \in I} \zeta_i e_i$ and $\nabla_\mathcal{D} \phi = \sum_{i \in I} \zeta_i \nabla e_i$. In what follows when we consider sequences $(\mathcal{D}_n)_{n \in \mathbb{N}}$ of gradient discretizations, it is useful to think that each $\mathcal{D}_n$ is associated to a mesh of size $h_n$ with $\lim_{n \to \infty} h_n = 0$.*

In the following $(\mathcal{D}_n)_{n \in \mathbb{N}}$ is a sequence of gradient discretizations.

**Definition 2.5.** *If $\mathcal{D}$ is a GD, define $C_\mathcal{D}$ as the norm of $\Pi_\mathcal{D}$:*

$$C_\mathcal{D} := \max_{\phi \in X_\mathcal{D} \setminus \{0\}} \frac{\|\Pi_\mathcal{D} \phi\|_{L^2(\Omega)}}{\|\nabla_\mathcal{D} \phi\|_{L^2(\Omega)^d}}.$$

*A sequence $(\mathcal{D}_n)_{n \in \mathbb{N}}$ of GD is coercive if there exists $C_p \in \mathbb{R}_+$ such that $C_{\mathcal{D}_n} \leq C_p$ for all $n \in \mathbb{N}$.*

We observe that coercivity implies a kind of Poincaré inequality.

**Definition 2.6.** *If $\mathcal{D}$ is a GD, define $S_\mathcal{D} : H^1_0(\Omega) \to [0, \infty[$ by*

$$S_\mathcal{D}(v) := \min_{\phi \in X_\mathcal{D}} (\|\Pi_\mathcal{D} \phi - v\|_{L^2(\Omega)} + \|\nabla_\mathcal{D} \phi - \nabla v\|_{L^2(\Omega)^d}).$$

*A sequence $(\mathcal{D}_n)_{n \in \mathbb{N}}$ of GD is consistent if $\lim_{n \to \infty} S_{\mathcal{D}_n}(v) = 0$ for all $v \in H^1_0(\Omega)$.*

**Definition 2.7.** *If $\mathcal{D}$ is a GD, define $W_\mathcal{D} : H_{\text{div}}(\Omega) \to [0, \infty[$ by*

$$W_\mathcal{D}(\boldsymbol{v}) = \sup_{\phi \in X_\mathcal{D} \setminus \{0\}} \frac{\left|\int_\Omega (\nabla_\mathcal{D} \phi \cdot \boldsymbol{v} + \Pi_\mathcal{D} \phi \nabla \cdot \boldsymbol{v}) d\boldsymbol{x}\right|}{\|\nabla_\mathcal{D} \phi\|_{L^2(\Omega)^d}}.$$

*A sequence $(\mathcal{D}_n)_{n \in \mathbb{N}}$ of GD is limit-conforming if $\lim_{n \to \infty} W_{\mathcal{D}_n}(\boldsymbol{v}) = 0$ for all $\boldsymbol{v} \in H_{\text{div}}(\Omega)$.*

The limit conformity of the method implies that the gradient discretization method satisfies asymptotically the Stokes theorem.

**Definition 2.8.** *A sequence $(\mathcal{D}_n)_{n \in \mathbb{N}}$ of GD is compact if, for any sequence $\phi_n \in X_{\mathcal{D}_n}$ such that $(\|\nabla_{\mathcal{D}_n} \phi_n\|_{L^2(\Omega)^d})_{n \in \mathbb{N}}$ is bounded, the sequence $(\Pi_{\mathcal{D}_n} \phi_n)_{n \in \mathbb{N}}$ is relatively compact in $L^2(\Omega)$.*

In order to use the GD to solve Eq. (2.1) it is useful to write $f \in H^{-1}(\Omega)$ as

$$f = f_0 + \sum_{i=1}^{d} \frac{\partial f_i}{\partial x_i} = f_0 + \nabla \cdot \boldsymbol{F},$$

where $\boldsymbol{x} = (x_1, ..., x_d) \in \Omega$, $f_0, f_1, ..., f_d \in L^2(\Omega)$ and $\boldsymbol{F} = (f_1, ..., f_d)^\top \in L^2(\Omega)^d$. With this notation, Eq. (2.1) becomes

$$\int_\Omega A\nabla u \cdot \nabla v \, d\boldsymbol{x} = \int_\Omega (f_0 \, v - \boldsymbol{F} \cdot \nabla v) \, d\boldsymbol{x} \qquad \text{for all } v \in H^1_0(\Omega). \tag{2.2}$$

We next define the Gradient Scheme used to approximate $u$ solution of Eq. (2.2).

**Definition 2.9.** *For a given gradient discretization $\mathcal{D}$, the Gradient Scheme (GS) for problem Eq. (2.2) is defined by: find $\vartheta \in X_\mathcal{D}$ such that*

$$\int_\Omega A\nabla_\mathcal{D}\vartheta \cdot \nabla_\mathcal{D}\phi \, d\boldsymbol{x} = \int_\Omega (f_0 \, \Pi_\mathcal{D}\phi - \boldsymbol{F} \cdot \nabla_\mathcal{D}\phi) \, d\boldsymbol{x} \qquad \text{for all } \phi \in X_\mathcal{D}. \tag{2.3}$$

The convergence of the above scheme is given by Theorem 2.10, which is proven in [9, Theorem 2.28]. Notice that under Assumption 2.1 and $u \in H^1_0(\Omega)$ we have $A\nabla u + \boldsymbol{F} \in H_{\text{div}}(\Omega)$, indeed $-\nabla \cdot (A\nabla u + \boldsymbol{F}) = f_0 \in L^2(\Omega)$.

**Theorem 2.10.** *Let $\mathcal{D}$ be a GD, then there exists one and only one $\vartheta \in X_\mathcal{D}$ solution to Eq. (2.3) and it satisfies*

$$\|\nabla u - \nabla_\mathcal{D}\vartheta\|_{L^2(\Omega)^d} \leq \frac{1}{\underline{\lambda}} W_\mathcal{D}(A\nabla u + \boldsymbol{F}) + (1 + \kappa(A))S_\mathcal{D}(u), \tag{2.4}$$

*where $\kappa(A) = \overline{\lambda}/\underline{\lambda}$ is the condition number of $A$.*

**Corollary 2.11.** *If $(\mathcal{D}_n)_{n\in\mathbb{N}}$ is a consistent and limit-conforming sequence of GD and $\vartheta_n \in \mathcal{D}_n$ is a sequence of solutions to Eq. (2.3), then*

$$\lim_{n\to\infty} \|\nabla u - \nabla_{\mathcal{D}_n}\vartheta_n\|_{L^2(\Omega)^d} = 0.$$

*Proof.* Follows from Eq. (2.4) and the definitions of consistency and limit conformity. □

Convergence rates are obtained under stronger regularity hypothesis on the data and the solution, upon the introduction of a mesh and depend on the underlying discretization method. We refer to Corollary 2.17 at the end of Section 2.2 for such results. The compactness hypothesis of Definition 2.8 is needed to establish convergence of the gradient scheme when applied to nonlinear problems.

## 2.2 The Symmetric Weighted Discontinuous Galerkin Gradient Discretizazion

Inspired from the method proposed in [11] we define the Symmetric Weighted Discontinuous Galerkin GD (SWDGGD).

A polytopal mesh $\mathfrak{T} = (\mathcal{M}, \mathcal{F}, \mathcal{P})$ is defined as follows. $\mathcal{M}$ is a finite family of non empty polytopal open disjoint elements $K \subset \Omega$ such that $\overline{\Omega} = \cup_{K \in \mathcal{M}} \overline{K}$. We suppose that $K$ is star shaped with respect to an $\boldsymbol{x}_K \in K$ and denote $\mathcal{P} = (\boldsymbol{x}_K)_{K \in \mathcal{M}}$. Let $\mathcal{F} = \mathcal{F}_b \cup \mathcal{F}_i$ be the set of faces of the mesh, where $\mathcal{F}_b$, $\mathcal{F}_i$ are the boundary and internal faces, respectively. The set of faces of $K$ is $\mathcal{F}_K = \{\sigma \in \mathcal{F} : \sigma \subset \partial K\}$. For each $K \in \mathcal{M}$ and $\sigma \in \mathcal{F}_K$ we denote by $d_{K,\sigma}$ the orthogonal distance between $\bm{x}_K$ and $\sigma$, hence

$$d_{K,\sigma} = (\bm{y} - \bm{x}_K) \cdot \bm{n}_{K,\sigma} \quad \text{for all } \bm{y} \in \sigma,$$

where $\bm{n}_{K,\sigma}$ is the unit vector normal to $\sigma$ outward to $K$. We denote by $D_{K,\sigma}$ the cone with vertex $\bm{x}_K$ and basis $\sigma$, that is

$$D_{K,\sigma} = \{\bm{x}_K + s(\bm{y} - \bm{x}_K) : s \in ]0,1[,\, \bm{y} \in \sigma\}.$$

Finally, we define the mesh size and a constant measuring the regularity of the mesh. For $\sigma \in \mathcal{F}$ let $\mathcal{M}_\sigma = \{K \in \mathcal{M} : \sigma \in \mathcal{F}_K\}$ and let $h_K$ be the diameter of $K \in \mathcal{M}$, then

$$h_\mathcal{M} = \max\{h_K : K \in \mathcal{M}\},$$
$$\eta_{\mathfrak{T}} = \max\left(\{\frac{h_T}{h_K} + \frac{h_K}{h_T} : \sigma \in \mathcal{F}_i, \mathcal{M}_\sigma = \{K,T\}\} \cup \{\frac{h_K}{d_{K,\sigma}} : K \in \mathcal{M}, \sigma \in \mathcal{F}_K\}\right.$$
$$\left.\cup \{\#\mathcal{F}_K : K \in \mathcal{M}\}\right),$$

the term $\{\#\mathcal{F}_K : K \in \mathcal{M}\}$ is needed in [11, Lemma 3.14] to bound the jumps on the faces of the elements.

Let $V = \{v \in L^2(\Omega) : v|_K \in \mathbb{P}_\ell(K), \forall K \in \mathcal{M}\}$, where $\mathbb{P}_\ell(K)$ is the space of polynomials in $K$ of total degree $\ell$. Let $(e_i)_{i \in I}$ be a basis of $V$ such that $\mathrm{supp}(e_i)$ is restricted to one element of $\mathcal{M}$. We set

$$X_\mathcal{D} = \{\phi = (\zeta_i)_{i \in I} : \zeta_i \in \mathbb{R} \text{ for all } i \in I\} \tag{2.5}$$

and define the operator $\Pi_\mathcal{D} : X_\mathcal{D} \to L^2(\Omega)$ by

$$\Pi_\mathcal{D}\phi = \sum_{i \in I} \zeta_i e_i. \tag{2.6}$$

For $K \in \mathcal{M}$ we note by $\Pi_{\overline{K}}\phi := \Pi_\mathcal{D}\phi|_{\overline{K}}$ the restriction of $\Pi_\mathcal{D}\phi$ to $K$ extended to $\overline{K}$ and define $\nabla_{\overline{K}}\phi = \nabla\Pi_{\overline{K}}\phi$. Let $\alpha \in ]0,1[$ be a user parameter and $\psi : [0,1] \to \mathbb{R}$ such that $\psi(s) = 0$ on $[0,\alpha[$ and $\psi|_{[\alpha,1]} \in \mathbb{P}_{\ell-1}([\alpha,1])$ satisfying

$$\int_\alpha^1 \psi(s) s^{d-1} ds = 1 \quad \text{and} \quad \int_\alpha^1 (1-s)^i \psi(s) s^{d-1} ds = 0 \quad \text{for } i = 1, ..., \ell - 1. \tag{2.7}$$

In the case where $\ell = 1$ we have $\psi(s)|_{[\alpha,1]} = d/(1-\alpha^d)$. This choice of $\psi$ is fundamental to show the equivalence with the SWIPG method, see Appendix A. The discrete gradient $\nabla_\mathcal{D} : X_\mathcal{D} \to L^2(\Omega)^d$ is defined as follows. For $\phi \in X_\mathcal{D}$, $K \in \mathcal{M}$ and $\sigma \in \mathcal{F}_K$, we set, for a.e. $\bm{x} \in D_{K,\sigma}$

$$\nabla_\mathcal{D}\phi(\bm{x}) = \nabla_{\overline{K}}\phi(\bm{x}) + \psi(s)\frac{[\phi]_{K,\sigma}(\bm{y})}{d_{K,\sigma}}\bm{n}_{K,\sigma}, \tag{2.8}$$

where $\bm{x} = \bm{x}_K + s(\bm{y} - \bm{x}_K)$ with $s \in ]0,1[$, $\bm{y} \in \sigma$ and

if $\sigma \in \mathcal{F}_i$ with $\sigma = \partial K \cap \partial T$ then $[\phi]_{K,\sigma}(\bm{y}) = \omega_{K,\sigma}(\Pi_{\overline{T}}\phi(\bm{y}) - \Pi_{\overline{K}}\phi(\bm{y}))$,
if $\sigma \in \mathcal{F}_b$ with $\sigma = \partial K \cap \partial \Omega$ then $[\phi]_{K,\sigma}(\bm{y}) = 0 - \Pi_{\overline{K}}\phi(\bm{y})$.

For $\sigma \in \mathcal{F}_b$ with $\sigma = \partial K \cap \partial \Omega$ and $K \in \mathcal{M}$ it is useful to set $\omega_{K,\sigma} = 1$. If instead $\sigma \in \mathcal{F}_i$ with $\sigma = \partial K \cap \partial T$ and $K, T \in \mathcal{M}$ the weights $\omega_{K,\sigma}, \omega_{T,\sigma}$ are two non negative numbers such that

$$\omega_{K,\sigma} + \omega_{T,\sigma} = 1. \tag{2.9}$$

In the original Discontinuous Galerkin GD (DGGD) method introduced in [11] the weights are $(\omega_{K,\sigma}, \omega_{T,\sigma}) = (1/2, 1/2)$ and it is proven that $\mathcal{D} = (X_\mathcal{D}, \Pi_\mathcal{D}, \nabla_\mathcal{D})$, with $X_\mathcal{D}, \Pi_\mathcal{D}, \nabla_\mathcal{D}$ as in Eqs. (2.5), (2.6) and (2.8), is a GD method. Moreover, any sequence $(\mathcal{D}_n)_{n\in\mathbb{N}}$ of DGGD defined from polytopal meshes $(\mathfrak{T}_n)_{n\in\mathbb{N}}$ with $(\eta_{\mathfrak{T}_n})_{n\in\mathbb{N}}$ bounded and $h_{\mathcal{M}_n} \longrightarrow 0$ is a coercive, consistent, limit-conforming and compact sequence of GD. Thanks to the particular choice of $\psi$ it is possible to show that in the linear case with piecewise constant diffusion the DGGD scheme is equivalent to the well known SIPG method.

In our case we want to be equivalent to the SWIPG method, hence we define the weights as follows. Let $K \in \mathcal{M}$ and $\sigma \in \mathcal{F}_K$, we set

$$\delta_{K,\sigma} = \boldsymbol{n}_{K,\sigma}^\top A|_K \boldsymbol{n}_{K,\sigma}.$$

For $\sigma \in \mathcal{F}_i$ such that $\sigma = \partial K \cap \partial T$ with $K, T \in \mathcal{M}$ we define

$$\omega_{K,\sigma} = \frac{\delta_{T,\sigma}}{\delta_{K,\sigma} + \delta_{T,\sigma}}, \qquad \omega_{T,\sigma} = \frac{\delta_{K,\sigma}}{\delta_{K,\sigma} + \delta_{T,\sigma}}. \tag{2.10}$$

Upon changing the constants in [11, Lemma 3.8] we deduce from [11, Lemma 3.10] that $\|\nabla_\mathcal{D} \cdot \|_{L^2(\Omega)^d}$ with the choice of weights given by Eq. (2.10) is a norm on $X_\mathcal{D}$ and hence $\mathcal{D} = (X_\mathcal{D}, \Pi_\mathcal{D}, \nabla_\mathcal{D})$ with $(\omega_{K,\sigma}, \omega_{T,\sigma})$ as in Eq. (2.10) is a GD method. It can be used to solve diffusion problems with homogeneous boundary conditions as in Definition 2.9. From now on we refer to this method as the Symmetric Weighted DGGD scheme (SWDGGD). Apart from the weights definition, the only difference with respect to the original DGGD method is a factor

$$C_\omega := \max_{K \in \mathcal{M}, \sigma \in \mathcal{F}_K} \omega_{K,\sigma}^{-1}/2 \tag{2.11}$$

multiplying the constant $C_\mathcal{D}$ of Definition 2.5.

In the foregoing analysis we need the jump semi norm on $X_\mathcal{D}$, defined by

$$|\phi|_J^2 := \sum_{K \in \mathcal{M}} \sum_{\sigma \in \mathcal{F}_K} \frac{1}{d_{K,\sigma}} \int_\sigma [\phi]_{K,\sigma}^2(\boldsymbol{y}) d\boldsymbol{y}.$$

We define a stronger version of $S_\mathcal{D}$ which controls the jumps.

**Definition 2.12.** If $\mathcal{D}$ is a SWDGGD, define $S_{\mathcal{D},J} : H_0^1(\Omega) \to [0, \infty[$ by

$$S_{\mathcal{D},J}(v) := \min_{\phi \in X_\mathcal{D}} (\|\Pi_\mathcal{D}\phi - v\|_{L^2(\Omega)} + \|\nabla_\mathcal{D}\phi - \nabla v\|_{L^2(\Omega)^d} + |\phi|_J).$$

We quote two improved estimates on $S_\mathcal{D}$, $S_{\mathcal{D},J}$ and $W_\mathcal{D}$.

**Lemma 2.13.** *There exists $C_S > 0$ depending only on $|\Omega|$, $\alpha$, $\ell$, $d$ and $\rho$ such that for all $v \in H^2(\Omega) \cap H_0^1(\Omega)$*

$$S_\mathcal{D}(v) \leq C_S h_\mathcal{M} \|v\|_{H^2(\Omega)} \quad \text{and} \quad S_{\mathcal{D},J}(v) \leq C_S h_\mathcal{M} \|v\|_{H^2(\Omega)}.$$

The result for $S_{\mathcal{D},J}$ is obtained following the lines of the proof for $S_\mathcal{D}$, which is given in [11, Lemma 3.14].

**Lemma 2.14.** *There exists $C_W > 0$ depending only on $|\Omega|$, $\alpha$, $\ell$ and $d$ such that for all $\boldsymbol{v} \in H^1(\Omega)^d$*

$$W_{\mathcal{D}}(\boldsymbol{v}) \leq C_W h_{\mathcal{M}} \|\boldsymbol{v}\|_{H^1(\Omega)^d}.$$

Lemma 2.14 has been proven for the DGGD scheme in [11, Lemma 3.15]. The proof uses the fact that $(1/2, 1/2)$ is a partition of unity. Thanks to Eq. (2.9) the same result holds for the SWDGGD method. Next, Theorem 2.15 establishes the asymptotic properties of the SWDGGD schemes.

**Theorem 2.15.** *Let $(\mathcal{D}_n)_{n \in \mathbb{N}}$ be a sequence of of SWDGGD defined from polytopal meshes $(\mathfrak{T}_n)_{n \in \mathbb{N}}$ with $(\eta_{\mathfrak{T}_n})_{n \in \mathbb{N}}$ bounded and $h_{\mathcal{M}_n} \longrightarrow 0$ for $n \to \infty$. Then it is a coercive, consistent, limit-conforming and compact sequence of GD.*

*Proof.* Coercivity and compactness are proven as in [11, Lemma 3.12, Lemma 3.13]. Consistency follows from Lemma 2.13 and [9, Lemma 2.16]. Limit-conformity follows from the compactness of the scheme, Lemma 2.14 and [9, Lemma 2.17]. □

In the SWDGGD scheme the $C_{\mathrm{p}}$ constant depends continuously on $C_\omega$ from Eq. (2.11). We note that, even if $C_\omega$ is mesh dependent it can be bounded by terms depending only on $A$. In the following Lemma we show, by usual density arguments, that even if $v$ is only in $H_0^1(\Omega)$ we have $\lim_{n \to \infty} S_{\mathcal{D}_n, J}(v) = 0$.

**Lemma 2.16.** *Consider the same assumptions of Theorem 2.15 and $v \in H_0^1(\Omega)$. Then we have $\lim_{n \to \infty} S_{\mathcal{D}_n, J}(v) = 0$.*

*Proof.* Let $v \in H_0^1(\Omega)$ and $\varepsilon > 0$. Then there exists $v_\varepsilon \in H^2(\Omega) \cap H_0^1(\Omega)$ such that $\|v - v_\varepsilon\|_{L^2(\Omega)} + \|\nabla v - \nabla v_\varepsilon\|_{L^2(\Omega)^d} \leq \varepsilon$. Let

$$\phi_n = \operatorname*{argmin}_{\phi \in X_{\mathcal{D}_n}} (\|\Pi_{\mathcal{D}_n} \phi - v_\varepsilon\|_{L^2(\Omega)} + \|\nabla_{\mathcal{D}_n} \phi - \nabla v_\varepsilon\|_{L^2(\Omega)^d} + |\phi|_J).$$

Hence

$$\begin{aligned} S_{\mathcal{D}_n, J}(v) \leq & \|\Pi_{\mathcal{D}_n} \phi_n - v\|_{L^2(\Omega)} + \|\nabla_{\mathcal{D}_n} \phi_n - \nabla v\|_{L^2(\Omega)^d} + |\phi_n|_J \\ \leq & \varepsilon + C_S h_{\mathcal{M}_n} \|v_\varepsilon\|_{H^2(\Omega)}, \end{aligned}$$

using Lemma 2.13 $\lim_{n \to \infty} S_{\mathcal{D}_n, J}(v) \leq \varepsilon$. Since $\varepsilon$ is arbitrary the result follows. □

**Corollary 2.17** (Of Theorem 2.10). *Let $\mathcal{D}$ be a SWDGGD, under the same assumptions of Theorem 2.10, $u \in H^2(\Omega) \cap H_0^1(\Omega)$ and $\boldsymbol{F} \in H^1(\Omega)^d$, the solution $\vartheta \in X_{\mathcal{D}}$ to Eq. (2.3) satisfies*

$$\|\nabla u - \nabla_{\mathcal{D}} \vartheta\|_{L^2(\Omega)^d} \leq h_{\mathcal{M}} \left( \frac{1}{\underline{\lambda}} C_W \|A \nabla u + \boldsymbol{F}\|_{H^1(\Omega)^d} + (1 + \kappa(A)) C_S \|u\|_{H^2(\Omega)} \right),$$

*Proof.* Follows from Theorem 2.10 together with Lemmas 2.13 and 2.14. □

## 2.3 The Local Weighted Discontinuous Galerkin Gradient Discretization

Let $\{\Omega_k\}_{k=1}^M$ be a sequence of polytopal domains with $\Omega_1 = \Omega$ and $\Omega_k \subset \Omega$. We consider as well a sequence $(\mathfrak{T}_k)_{k=1}^M = ((\mathcal{M}_k, \mathcal{F}_k, \mathcal{P}_k))_{k=1}^M$ of polytopal meshes on $\Omega$ and denote $\mathcal{F}_k = \mathcal{F}_{k,b} \cup \mathcal{F}_{k,i}$ with $\mathcal{F}_{k,b}$ and $\mathcal{F}_{k,i}$ the set of boundary and internal faces of $\mathcal{M}_k$. Moreover, $(\mathfrak{T}_k)_{k=1}^M$ satisfies the following.

**Assumption 2.18.**

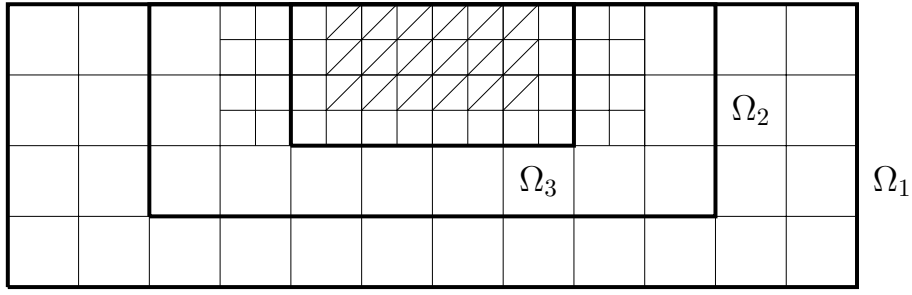

Figure 1. *Example of possible meshes for three embedded domains $\Omega_1$, $\Omega_2$, $\Omega_3$.*

a) For each $k = 1, ..., M$, $\overline{\Omega}_k = \cup_{K \in \mathcal{M}_k,\, K \subset \Omega_k} \overline{K}$.

b) For $k = 1, ..., M - 1$

   i) $\{K \in \mathcal{M}_{k+1} : K \subset \Omega \setminus \Omega_{k+1}\} = \{K \in \mathcal{M}_k : K \subset \Omega \setminus \Omega_{k+1}\}$,

   ii) if $K, T \in \mathcal{M}_k$ with $K \subset \Omega_{k+1}$, $T \subset \Omega \setminus \Omega_{k+1}$ and $\partial K \cap \partial T \neq \emptyset$ then $K \in \mathcal{M}_{k+1}$,

   iii) if $K \in \mathcal{M}_k$ and $K \subset \Omega_{k+1}$, either $K \in \mathcal{M}_{k+1}$ or $K$ is a union of elements in $\mathcal{M}_{k+1}$.

c) We suppose the existence of $C_r > 0$ such that

   i) for $k = 1, ..., M - 1$, if $K \in \mathcal{M}_k$ and $\widehat{K} \in \mathcal{M}_{k+1}$ with $\widehat{K} \subset K$ and $\sigma \in \mathcal{F}_K$, $\hat{\sigma} \in \mathcal{F}_{\widehat{K}}$ with $\hat{\sigma} \subset \sigma$ then $d_{K,\sigma} \leq C_r d_{\widehat{K},\hat{\sigma}}$,

   ii) for $k = 1, ..., M$, if $\sigma = \partial K \cap \partial T$ with $K, T \in \mathcal{M}_k$, $T \subset \Omega \setminus \Omega_k$ and $K \subset \Omega_k$ then $d_{K,\sigma} \leq C_r d_{T,\sigma}$.

d) It exists $\rho > 0$ such that $\eta_{\mathfrak{T}_k} \leq \rho$ for $k = 1, ..., M$.

The above assumptions on $(\mathfrak{T}_k)_{k=1}^M$ ensure that $\mathfrak{T}_{k+1}$ is a refinement of $\mathfrak{T}_k$ and that this refinement occurs inside the subdomain $\Omega_{k+1}$. Let $\widehat{\mathfrak{T}}_k = (\widehat{\mathcal{M}}_k, \widehat{\mathcal{F}}_k, \widehat{\mathcal{P}}_k)$, with $\widehat{\mathcal{M}}_k = \{K \in \mathcal{M}_k : K \subset \Omega_k\}$, $\widehat{\mathcal{P}}_k = \{\boldsymbol{x}_k \in \mathcal{P}_k : \boldsymbol{x}_k \in \Omega_k\}$ and $\widehat{\mathcal{F}}_k = \widehat{\mathcal{F}}_{k,b} \cup \widehat{\mathcal{F}}_{k,i}$ the set of faces of $\widehat{\mathcal{M}}_k$, with $\widehat{\mathcal{F}}_{k,b}$ and $\widehat{\mathcal{F}}_{k,i}$ the boundary and internal faces of $\Omega_k$, respectively. Condition a) in Assumption 2.18 assures that $\widehat{\mathfrak{T}}_k$ is a polytopal mesh on $\Omega_k$. b) guarantees that in $\Omega \setminus \Omega_{k+1}$ and in the neighborhood of $\partial \Omega_{k+1}$ the meshes $\mathcal{M}_k$ and $\mathcal{M}_{k+1}$ are equal and that $\mathcal{M}_{k+1}$ is a refinement of $\mathcal{M}_k$ in $\Omega_{k+1}$. c) and d) ensure mesh regularity, will permit equivalences between jump norms and make the constant $C_S$ of Lemma 2.13 uniform in $k$. An example of meshes satisfying Assumption 2.18 is given in Section 2.3.

Given $(\mathfrak{T}_k)_{k=1}^M$ we define a sequence $\mathcal{D}_k = (X_{\mathcal{D}_k}, \Pi_{\mathcal{D}_k}, \nabla_{\mathcal{D}_k})$ of SWDGGD. Let

$$V_k = \{v_k \in L^2(\Omega) : v_k|_K \in \mathbb{P}_\ell(K), \forall K \in \mathcal{M}_k\} \tag{2.12}$$

and $(e_{k,i})_{i \in I_k}$ be a basis of $V_k$ such that $\mathrm{supp}(e_{k,i})$ is restricted to one element of $\mathcal{M}_k$. We set

$$X_{\mathcal{D}_k} = \{\phi_k = (\zeta_{k,i})_{i \in I_k} : \zeta_{k,i} \in \mathbb{R} \text{ for all } i \in I_k\}.$$

$\Pi_{\mathcal{D}_k}$ and $\nabla_{\mathcal{D}_k}$ are defined as in Eqs. (2.6), (2.8) and (2.10).

We can write $X_{\mathcal{D}_k} = Y_{\mathcal{D}_k} \oplus Z_{\mathcal{D}_k}$, where $\text{supp}(\Pi_{\mathcal{D}_k}\varphi_k) \subset \Omega_k$ for $\varphi_k \in Y_{\mathcal{D}_k}$ and $\text{supp}(\Pi_{\mathcal{D}_k}\xi_k) \subset \Omega \setminus \Omega_k$ for $\xi_k \in Z_{\mathcal{D}_k}$. For $k = 1$ we have $Y_{\mathcal{D}_1} = X_{\mathcal{D}_1}$ and $Z_{\mathcal{D}_1} = \{0\}$. For $k \geq 2$ and $\phi_{k-1} \in X_{\mathcal{D}_{k-1}}$ there exists $\xi_k \in Z_{\mathcal{D}_k}$ such that $\Pi_{\mathcal{D}_{k-1}}\phi_{k-1}\chi_{\Omega\setminus\Omega_k} = \Pi_{\mathcal{D}_k}\xi_k$. By abuse of notation we we will denote $\xi_k = \phi_{k-1}\chi_{\Omega\setminus\Omega_k}$, hence $\chi_{\Omega\setminus\Omega_k}$ is seen as an operator from $X_{\mathcal{D}_{k-1}}$ to $Z_{\mathcal{D}_k}$.

In what follows $\Pi_{\widehat{\mathcal{D}}_k}$ is the restriction of $\Pi_{\mathcal{D}_k}$ to $Y_{\mathcal{D}_k}$. Let us define as well a gradient on $Y_{\mathcal{D}_k}$ which will be used to impose inhomogeneous Dirichlet boundary conditions. Let $\varphi_k \in Y_{\mathcal{D}_k}$ and $\xi_k \in Z_{\mathcal{D}_k}$, for $K \in \widehat{\mathcal{M}}_k$, $\sigma \in \mathcal{F}_K$ and $\boldsymbol{x} \in D_{K,\sigma}$ the gradient $\nabla_{\widehat{\mathcal{D}}_k,\xi_k}\varphi_k(\boldsymbol{x})$ is defined by

$$\nabla_{\widehat{\mathcal{D}}_k,\xi_k}\varphi_k(\boldsymbol{x}) = \nabla_{\overline{K}}\varphi_k(\boldsymbol{x}) + \psi(s)\frac{[\varphi_k]_{K,\sigma,\xi_k}(\boldsymbol{y})}{d_{K,\sigma}}\boldsymbol{n}_{K,\sigma},$$

where $\boldsymbol{x} = \boldsymbol{x}_K + s(\boldsymbol{y} - \boldsymbol{x}_K)$ and

$$[\varphi_k]_{K,\sigma,\xi_k}(\boldsymbol{y}) = [\varphi_k]_{K,\sigma}(\boldsymbol{y}) \quad \text{if } \sigma \in \widehat{\mathcal{F}}_{k,i} \text{ or } \sigma \in \widehat{\mathcal{F}}_{k,b} \cap \mathcal{F}_{k,b},$$

$$[\varphi_k]_{K,\sigma,\xi_k}(\boldsymbol{y}) = \Pi_{\overline{T}}\xi_k - \Pi_{\overline{K}}\varphi_k \quad \text{if } \sigma \in \widehat{\mathcal{F}}_{k,b} \setminus \mathcal{F}_{k,b} \text{ with } \sigma = \partial K \cap \partial T$$
$$\text{and } K \in \widehat{\mathcal{M}}_k, T \in \mathcal{M}_k \setminus \widehat{\mathcal{M}}_k.$$

We will denote $\nabla_{\widehat{\mathcal{D}}_k,0}$ by $\nabla_{\widehat{\mathcal{D}}_k}$.

**Theorem 2.19.** *The triple $\widehat{\mathcal{D}}_k = (Y_{\mathcal{D}_k}, \Pi_{\widehat{\mathcal{D}}_k}, \nabla_{\widehat{\mathcal{D}}_k})$ is a SWDGGD scheme for each $k = 1, ..., M$.*

*Proof.* We notice that $\widehat{\mathcal{D}}_k$ is the SWDGGD corresponding to the local polytopal mesh $\widehat{\mathfrak{T}}_k$, hence it is a SWDGGD by construction. $\square$

In what follows we will call $\widehat{\mathcal{D}}_k$ the local SWDGGD. Remark that Lemmas 2.13 and 2.14 and Theorem 2.15 are valid if we replace $\mathcal{D}$, $\Omega$, $h_\mathcal{M}$ and $\mathfrak{T}$ with $\widehat{\mathcal{D}}_k$, $\Omega_k$, $h_{\widehat{\mathcal{M}}_k}$ and $\widehat{\mathfrak{T}}_k$.

Observe that for $\varphi_k \in Y_{\mathcal{D}_k}$ we have $\nabla_{\widehat{\mathcal{D}}_k}\varphi_k \neq \nabla_{\mathcal{D}_k}\varphi_k$, indeed $\nabla_{\widehat{\mathcal{D}}_k}$ is missing the $\omega_{K,\sigma}$ factor in the jumps at the faces $\sigma \in \widehat{\mathcal{F}}_{k,b} \setminus \mathcal{F}_{k,b}$. Adding the $\omega_{K,\sigma}$ factor in those faces would prevent the limit consistency of $\widehat{\mathcal{D}}_k$.

In what follows $S_{\widehat{\mathcal{D}}_k}$ and $W_{\widehat{\mathcal{D}}_k}$ are the operators defined by Definitions 2.6 and 2.7 but with $\Omega$, $\mathcal{D}$, and $X_\mathcal{D}$ replaced by $\Omega_k$, $\widehat{\mathcal{D}}_k$ and $Y_{\mathcal{D}_k}$. We define as well the jump semi norms on $X_{\mathcal{D}_k}$ and $Y_{\mathcal{D}_k}$. For $\phi_k \in X_{\mathcal{D}_k}$ we define

$$|\phi_k|^2_{J(k)} := \sum_{K \in \mathcal{M}_k} \sum_{\sigma \in \mathcal{F}_K} \frac{1}{d_{K,\sigma}} \int_\sigma [\phi_k]_{K,\sigma}(\boldsymbol{y})^2 d\boldsymbol{y}$$

and for $\xi_k \in Z_{\mathcal{D}_k}$, $\varphi_k \in Y_{\mathcal{D}_k}$ we set

$$|\varphi_k|^2_{\widehat{J}(k),\xi_k} := \sum_{K \in \widehat{\mathcal{M}}_k} \sum_{\sigma \in \mathcal{F}_K} \frac{1}{d_{K,\sigma}} \int_\sigma [\varphi_k]_{K,\sigma,\xi_k}(\boldsymbol{y})^2 d\boldsymbol{y}.$$

Since in our local scheme (to be defined in Section 3) we solve local elliptic problems with artificial boundary conditions we need a local version of $S_{\mathcal{D}_k,J}$ which measures the error of the method on the boundary.

**Definition 2.20.** *Let $\xi_k \in Z_{\mathcal{D}_k}$ and $\widehat{\mathcal{D}}_k$ be a local SWDGGD, define $S_{\widehat{\mathcal{D}}_k,J,\xi_k} : H^1_0(\Omega) \to [0,\infty[$ by*

$$S_{\widehat{\mathcal{D}}_k,J,\xi_k}(v) := \min_{\varphi \in Y_{\mathcal{D}_k}} (\|\nabla_{\widehat{\mathcal{D}}_k,\xi_k}\varphi - \nabla v\|_{L^2(\Omega_k)^d} + |\varphi|_{\widehat{J}(k),\xi_k}).$$

The $L^2(\Omega_k)$ norm is not taken into account in $S_{\widehat{\mathcal{D}}_k,J,\xi_k}$ since our convergence results are in energy and jump norms.

**Lemma 2.21.** *Let $v \in H_0^1(\Omega) \cap H^2(\Omega)$, then for $k = 1, ..., M$*

$$\min_{\xi \in Z_{\mathcal{D}_k}} S_{\widehat{\mathcal{D}}_k,J,\xi}(v) \leq C_S h_{\widehat{\mathcal{M}}_k} \|v\|_{H^2(\Omega_k)}.$$

*Proof.* Follows the lines of [11, Lemma 3.14]. □

In order to provide bounds on $S_{\widehat{\mathcal{D}}_k,J,\xi_k}$ we need an additional norm to measure the error at the interface between subdomains. Let $\phi_k \in X_{\mathcal{D}_k}$, we define

$$|\phi_k|^2_{\partial\Omega_k^-} := \sum_{\substack{K\in\widehat{\mathcal{M}}_k \\ T\in\mathcal{M}_k\setminus\widehat{\mathcal{M}}_k}} \sum_{\sigma\in\mathcal{F}_K\cap\mathcal{F}_T} \frac{1}{d_{K,\sigma}} \int_\sigma \Pi_{\overline{T}}\phi_k(\boldsymbol{y})^2 \, d\boldsymbol{y}.$$

The minus in $|\cdot|_{\partial\Omega_k^-}$ refers to the fact that in the integral the argument lives outside $\widehat{\mathcal{M}}_k$. Later, $|\cdot|_{\partial\Omega_k^+}$ will be defined as well.

**Lemma 2.22.** *Let $\kappa_k, \xi_k \in Z_{\mathcal{D}_k}$ and $v \in H_0^1(\Omega)$. Then*

$$S_{\widehat{\mathcal{D}}_k,J,\kappa_k}(v) \leq S_{\widehat{\mathcal{D}}_k,J,\xi_k}(v) + C_\partial |\kappa_k - \xi_k|_{\partial\Omega_k^-},$$

*where $C_\partial = 1 + C_\psi$ and $C_\psi^2 = \int_\alpha^1 \psi(s)^2 s^{d-1} ds$. If moreover $v \in H^2(\Omega) \cap H_0^1(\Omega)$ we have*

$$S_{\widehat{\mathcal{D}}_k,J,\kappa_k}(v) \leq C_S h_{\widehat{\mathcal{M}}_k} \|v\|_{H^2(\Omega_k)} + C_\partial |\kappa_k - \xi_k|_{\partial\Omega_k^-} \quad \text{for} \quad \xi_k = \operatorname*{argmin}_{\xi \in Z_{\mathcal{D}_k}} S_{\widehat{\mathcal{D}}_k,J,\xi}(v).$$

*Proof.* Let $\kappa_k, \xi_k \in Z_{\mathcal{D}_k}$, $v \in H_0^1(\Omega)$ and $\varphi_k \in Y_{\mathcal{D}_k}$ defined by

$$\varphi_k = \operatorname*{argmin}_{\varphi \in Y_{\mathcal{D}_k}}(\|\nabla_{\widehat{\mathcal{D}}_k,\xi_k}\varphi - \nabla v\|_{L^2(\Omega_k)^d} + |\varphi|_{\widehat{J}(k),\xi_k}),$$

we have

$$\|\nabla_{\widehat{\mathcal{D}}_k,\kappa_k}\varphi_k - \nabla_{\widehat{\mathcal{D}}_k,\xi_k}\varphi_k\|^2_{L^2(\Omega_k)^d} = \sum_{K\in\widehat{\mathcal{M}}_k} \sum_{\sigma\in\mathcal{F}_K\cap\widehat{\mathcal{F}}_{k,b}} \int_{D_{K,\sigma}} \frac{\psi(s)^2}{d_{K,\sigma}^2} ([\varphi_k]_{K,\sigma,\kappa_k}(\boldsymbol{y}) - [\varphi_k]_{K,\sigma,\xi_k}(\boldsymbol{y}))^2 \, d\boldsymbol{x},$$

where $\boldsymbol{x} = \boldsymbol{x}_K + s(\boldsymbol{y} - \boldsymbol{x}_K)$ for $s \in [0,1]$ and $\boldsymbol{y} \in \sigma$. Using the change of variables $d\boldsymbol{x} = d_{K,\sigma} s^{d-1} ds d\boldsymbol{y}$ yields

$$\|\nabla_{\widehat{\mathcal{D}}_k,\kappa_k}\varphi_k - \nabla_{\widehat{\mathcal{D}}_k,\xi_k}\varphi_k\|^2_{L^2(\Omega_k)^d}$$
$$= \sum_{K\in\widehat{\mathcal{M}}_k} \sum_{\sigma\in\mathcal{F}_K\cap\widehat{\mathcal{F}}_{k,b}} \int_\sigma \int_\alpha^1 \frac{\psi(s)^2}{d_{K,\sigma}^2} ([\varphi_k]_{K,\sigma,\kappa_k}(\boldsymbol{y}) - [\varphi_k]_{K,\sigma,\xi_k}(\boldsymbol{y}))^2 d_{K,\sigma} s^{d-1} \, ds d\boldsymbol{y}$$
$$= C_\psi^2 \sum_{K\in\widehat{\mathcal{M}}_k} \sum_{\sigma\in\mathcal{F}_K\cap\widehat{\mathcal{F}}_{k,b}} \frac{1}{d_{K,\sigma}} \int_\sigma ([\varphi_k]_{K,\sigma,\kappa_k}(\boldsymbol{y}) - [\varphi_k]_{K,\sigma,\xi_k}(\boldsymbol{y}))^2 \, d\boldsymbol{y}.$$

If in the above sum $\sigma \in \mathcal{F}_{k,b}$, then $[\varphi_k]_{K,\sigma,\kappa_k} - [\varphi_k]_{K,\sigma,\xi_k} = 0$. Else, if $\sigma \in \mathcal{F}_{k,i}$ there is $T \in \mathcal{M}_k \setminus \widehat{\mathcal{M}}_k$ such that $\sigma = \partial K \cap \partial T$ and

$$[\varphi_k]_{K,\sigma,\kappa_k} - [\varphi_k]_{K,\sigma,\xi_k} = \Pi_{\overline{T}}\kappa_k - \Pi_{\overline{T}}\xi_k,$$

which implies

$$\|\nabla_{\widehat{\mathcal{D}}_k,\kappa_k}\varphi_k - \nabla_{\widehat{\mathcal{D}}_k,\xi_k}\varphi_k\|^2_{L^2(\Omega_k)^d}$$
$$= C_\psi^2 \sum_{\substack{K \in \widehat{\mathcal{M}}_k \\ T \in \mathcal{M}_k \setminus \widehat{\mathcal{M}}_k}} \sum_{\sigma \in \mathcal{F}_K \cap \mathcal{F}_T} \frac{1}{d_{K,\sigma}} \int_\sigma \Pi_{\overline{T}}(\kappa_k - \xi_k)(\boldsymbol{y})^2 \, d\boldsymbol{y} = C_\psi^2 |\kappa_k - \xi_k|^2_{\partial\Omega_k^-}. \qquad (2.13)$$

For the jump term $|\varphi_k|_{\widehat{J}(k),\kappa_k}$, we have

$$|\varphi_k|^2_{\widehat{J}(k),\kappa_k} = |\varphi_k|^2_{\widehat{J}(k),\xi_k}$$
$$+ \sum_{K \in \widehat{\mathcal{M}}_k} \sum_{\sigma \in \mathcal{F}_K \cap \widehat{\mathcal{F}}_{k,b}} \frac{1}{d_{K,\sigma}} \int_\sigma ([\varphi_k]_{K,\sigma,\kappa_k}(\boldsymbol{y})^2 - [\varphi_k]_{K,\sigma,\xi_k}(\boldsymbol{y})^2) \, d\boldsymbol{y}. \qquad (2.14)$$

If $\sigma \in \mathcal{F}_{k,b}$ then $[\varphi_k]^2_{K,\sigma,\kappa_k} - [\varphi_k]^2_{K,\sigma,\xi_k} = 0$, else, if $\sigma \in \mathcal{F}_{k,i}$ with $\sigma = \partial T \cap \partial K$, $T \in \mathcal{M}_k \setminus \widehat{\mathcal{M}}_k$ we have

$$[\varphi_k]^2_{K,\sigma,\kappa_k} - [\varphi_k]^2_{K,\sigma,\xi_k} = ([\varphi_k]_{K,\sigma,\kappa_k} - [\varphi_k]_{K,\sigma,\xi_k})([\varphi_k]_{K,\sigma,\kappa_k} + [\varphi_k]_{K,\sigma,\xi_k})$$
$$= (\Pi_{\overline{T}}\kappa_k - \Pi_{\overline{T}}\xi_k)(\Pi_{\overline{T}}\kappa_k - \Pi_{\overline{T}}\xi_k + 2[\varphi_k]_{K,\sigma,\xi_k}). \qquad (2.15)$$

Using Eqs. (2.14) and (2.15) we obtain

$$|\varphi_k|^2_{\widehat{J}(k),\kappa_k} \leq |\varphi_k|^2_{\widehat{J}(k),\xi_k} + \sum_{\substack{K \in \widehat{\mathcal{M}}_k \\ T \in \mathcal{M}_k \setminus \widehat{\mathcal{M}}_k}} \sum_{\sigma \in \mathcal{F}_K \cap \mathcal{F}_T} \frac{1}{d_{K,\sigma}} \int_\sigma \Pi_{\overline{T}}(\kappa_k - \xi_k)(\boldsymbol{y})^2 \, d\boldsymbol{y}$$
$$+ 2 \sum_{\substack{K \in \widehat{\mathcal{M}}_k \\ T \in \mathcal{M}_k \setminus \widehat{\mathcal{M}}_k}} \sum_{\sigma \in \mathcal{F}_K \cap \mathcal{F}_T} \frac{1}{d_{K,\sigma}} \int_\sigma |[\varphi_k]_{K,\sigma,\xi_k}(\boldsymbol{y})\Pi_{\overline{T}}(\kappa_k - \xi_k)(\boldsymbol{y})| \, d\boldsymbol{y}$$
$$\leq |\varphi_k|^2_{\widehat{J}(k),\xi_k} + |\kappa_k - \xi_k|^2_{\partial\Omega_k^-} + 2|\varphi_k|_{\widehat{J}(k),\xi_k}|\kappa_k - \xi_k|_{\partial\Omega_k^-}$$
$$= (|\varphi_k|_{\widehat{J}(k),\xi_k} + |\kappa_k - \xi_k|_{\partial\Omega_k^-})^2.$$

Using the above estimate and Eq. (2.13) we get

$$S_{\widehat{\mathcal{D}}_k,J,\kappa_k}(v) \leq \|\nabla_{\widehat{\mathcal{D}}_k,\kappa_k}\varphi - \nabla v\|_{L^2(\Omega_k)^d} + |\varphi_k|_{\widehat{J}(k),\kappa_k}$$
$$\leq \|\nabla_{\widehat{\mathcal{D}}_k,\xi_k}\varphi_k - \nabla v\|_{L^2(\Omega_k)^d} + |\varphi_k|_{\widehat{J}(k),\xi_k}$$
$$+ (1 + C_\psi)|\kappa_k - \xi_k|_{\partial\Omega_k^-}$$
$$= S_{\widehat{\mathcal{D}}_k,J,\xi_k}(v) + (1 + C_\psi)|\kappa_k - \xi_k|_{\partial\Omega_k^-}.$$

If moreover $v \in H^1_0(\Omega) \cap H^2(\Omega)$ and $\xi_k = \operatorname{argmin}_{\xi \in Z_{\mathcal{D}_k}} S_{\widehat{\mathcal{D}}_k,J,\xi}(v)$, Lemma 2.21 yields $S_{\widehat{\mathcal{D}}_k,J,\xi_k}(v) \leq C_S h_{\widehat{\mathcal{M}}_k}\|v\|_{H^2(\Omega_k)}$. □

# 3 The local elliptic scheme

We introduce here our local SWDGGD elliptic scheme before embarking into its a priori error analysis.

Set $\vartheta_0 = 0$ and define iteratively $\vartheta_k \in X_{\mathcal{D}_k}$ for $k = 1, ..., M$ as

$$\vartheta_k = \hat{\vartheta}_k + \kappa_k, \tag{3.1a}$$

where $\kappa_k \in Z_{\mathcal{D}_k}$ is defined as

$$\kappa_k = \vartheta_{k-1} \chi_{\Omega \setminus \Omega_k} \tag{3.1b}$$

and $\hat{\vartheta}_k \in Y_{\mathcal{D}_k}$ is the solution of the local problem

$$\int_{\Omega_k} A \nabla_{\widehat{\mathcal{D}}_k, \kappa_k} \hat{\vartheta}_k \cdot \nabla_{\widehat{\mathcal{D}}_k} \varphi \, d\boldsymbol{x} = \int_{\Omega_k} (f_0 \Pi_{\widehat{\mathcal{D}}_k} \varphi - \boldsymbol{F} \cdot \nabla_{\widehat{\mathcal{D}}_k} \varphi) \, d\boldsymbol{x} \tag{3.1c}$$

for all $\varphi \in Y_{\mathcal{D}_k}$.

Remember that $\nabla_{\widehat{\mathcal{D}}_k} = \nabla_{\widehat{\mathcal{D}}_k, 0}$, hence we use homogeneous boundary conditions for $\varphi$. Due to the definition of $\nabla_{\widehat{\mathcal{D}}_k, \kappa_k}$ the inhomogeneous Dirichlet boundary condition $\kappa_k$ is weakly imposed on $\hat{\vartheta}_k$. We have $\kappa_1 = 0$, hence $\vartheta_1 = \hat{\vartheta}_1 \in X_{\mathcal{D}_1}$. Then, for $k \geq 2$ the scheme Eq. (3.1) computes a new local solution $\hat{\vartheta}_k$ on a refined mesh $\widehat{\mathcal{M}}_k$, where the boundary condition is inherited from the previous solution $\vartheta_{k-1}$.

In Section 3.1 we perform the a priori error analysis for the local solutions $\hat{\vartheta}_k$ and provide bounds for the errors in the local domains $\Omega_k$. Section 3.2 improves the results of Section 3.1 in a particular case, showing that the error due to artificial boundary conditions is of higher order. Finally, Section 3.3 provides error bounds for the global solution $\vartheta_k$.

## 3.1 A priori error analysis for the local solution

In this section we proceed with the a priori analysis of the local elliptic scheme presented in Section 3. Before proving convergence of the scheme we need the following interpolation result.

**Lemma 3.1.** Let $\xi_{k-1} \in Z_{\mathcal{D}_{k-1}}$, $\varphi_{k-1} \in Y_{\mathcal{D}_{k-1}}$ and $\xi_k = (\xi_{k-1} + \varphi_{k-1})\chi_{\Omega \setminus \Omega_k} \in Z_{\mathcal{D}_k}$. Then there exists $\varphi_k \in Y_{\mathcal{D}_k}$ such that

$$\|\nabla_{\widehat{\mathcal{D}}_k, \xi_k} \varphi_k - \nabla_{\widehat{\mathcal{D}}_{k-1}, \xi_{k-1}} \varphi_{k-1}\|_{L^2(\Omega_k)^d} \leq C_i |\varphi_{k-1}|_{\widehat{\mathcal{J}}(k-1), \xi_{k-1}}, \tag{3.2a}$$

$$|\varphi_k|_{\widehat{\mathcal{J}}(k), \xi_k} \leq C_i |\varphi_{k-1}|_{\widehat{\mathcal{J}}(k-1), \xi_{k-1}}, \tag{3.2b}$$

with $C_i = \sqrt{2} C_\psi (1 + C_{\omega,k}^2 C_r)^{1/2}$, $C_{\omega,k} = \max_{K \in \mathcal{M}_k, \sigma \in \mathcal{F}_K} \omega_{K,\sigma}^{-1}$, $C_\psi$ from Lemma 2.22 and $C_r$ from Assumption 2.18.

*Proof.* Since $\widehat{\mathcal{M}}_k$ is a refinement of $\widehat{\mathcal{M}}_{k-1}$ in $\Omega_k$, there exists $\varphi_k \in Y_{\mathcal{D}_k}$ such that $\Pi_{\widehat{\mathcal{D}}_k} \varphi_k = \Pi_{\widehat{\mathcal{D}}_{k-1}} \varphi_{k-1}|_{\Omega_k}$. Hence

$$\|\nabla_{\widehat{\mathcal{D}}_k, \xi_k} \varphi_k - \nabla_{\widehat{\mathcal{D}}_{k-1}, \xi_{k-1}} \varphi_{k-1}\|^2_{L^2(\Omega_k)^d}$$

$$\leq 2 \sum_{K \in \widehat{\mathcal{M}}_k} \sum_{\sigma \in \mathcal{F}_K} \int_{D_{K,\sigma}} \left| \psi(s) \frac{[\varphi_k]_{K,\sigma,\xi_k}(\boldsymbol{y})}{d_{K,\sigma}} \right|^2 d\boldsymbol{x}$$

$$+ 2 \sum_{\substack{K \in \mathcal{M}_{k-1} \\ K \subset \Omega_k}} \sum_{\sigma \in \mathcal{F}_K} \int_{D_{K,\sigma}} \left| \psi(s) \frac{[\varphi_{k-1}]_{K,\sigma,\xi_{k-1}}(\boldsymbol{y})}{d_{K,\sigma}} \right|^2 d\boldsymbol{x},$$

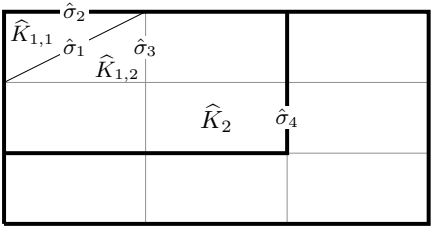

$\Omega_k$ and $\Omega_{k-1}$ bounded by thick lines ━,
$\hat{\sigma}_1 \subset K_1 = \widehat{K}_{1,1} \cup \widehat{K}_{1,2}$,
$\hat{\sigma}_2 = \sigma_2$ with $(\hat{\sigma}_2, \sigma_2) \in \widehat{\mathcal{F}}_{k,b} \times \widehat{\mathcal{F}}_{k-1,b}$,
$\hat{\sigma}_3 \in \widehat{\mathcal{F}}_{k,i}$,
$\hat{\sigma}_4 = \sigma_4$ with $(\hat{\sigma}_4, \sigma_4) \in \widehat{\mathcal{F}}_{k,b} \times \widehat{\mathcal{F}}_{k-1,i}$.

Figure 2. *Example of a situation described in the proof of Lemma 3.1.*

since the broken gradients of $\Pi_{\widehat{\mathcal{D}}_{k-1}}\varphi_{k-1}|_{\Omega_k}$ and $\Pi_{\widehat{\mathcal{D}}_k}\varphi_k$ cancel each other out. With the change of variables $d\boldsymbol{x} = d_{K,\sigma} s^{d-1} ds d\boldsymbol{y}$ we have

$$\int_{D_{K,\sigma}} \left| \psi(s) \frac{[\varphi_k]_{K,\sigma,\xi_k}(\boldsymbol{y})}{d_{K,\sigma}} \right|^2 d\boldsymbol{x} = \frac{1}{d_{K,\sigma}} \int_\alpha^1 \psi(s)^2 s^{d-1} ds \int_\sigma [\varphi_k]_{K,\sigma,\xi_k}(\boldsymbol{y})^2 d\boldsymbol{y}$$

and similarly for $\varphi_{k-1}$. Using $C_\psi^2 = \int_\alpha^1 \psi(s)^2 s^{d-1} ds$ yields

$$\|\nabla_{\widehat{\mathcal{D}}_k,\xi_k}\varphi_k - \nabla_{\widehat{\mathcal{D}}_{k-1},\xi_{k-1}}\varphi_{k-1}\|_{L^2(\Omega_k)^d}^2$$

$$\leq 2C_\psi^2 \sum_{K \in \widehat{\mathcal{M}}_k} \sum_{\sigma \in \mathcal{F}_K} \frac{1}{d_{K,\sigma}} \int_\sigma [\varphi_k]_{K,\sigma,\xi_k}(\boldsymbol{y})^2 d\boldsymbol{y}$$

$$+ 2C_\psi^2 \sum_{\substack{K \in \mathcal{M}_{k-1} \\ K \subset \Omega_k}} \sum_{\sigma \in \mathcal{F}_K} \frac{1}{d_{K,\sigma}} \int_\sigma [\varphi_{k-1}]_{K,\sigma,\xi_{k-1}}(\boldsymbol{y})^2 d\boldsymbol{y}$$

$$\leq 2C_\psi^2 (|\varphi_k|_{\widehat{J}(k),\xi_k}^2 + |\varphi_{k-1}|_{\widehat{J}(k-1),\xi_{k-1}}^2).$$

To obtain Eq. (3.2a) it remains to prove $|\varphi_k|_{\widehat{J}(k),\xi_k}^2 \leq C_{\omega,k}^2 C_{\mathrm{r}} |\varphi_{k-1}|_{\widehat{J}(k-1),\xi_{k-1}}^2$. We write $|\varphi_k|_{\widehat{J}(k),\xi_k}^2$ as

$$|\varphi_k|_{\widehat{J}(k),\xi_k}^2 = \sum_{\substack{K \in \mathcal{M}_{k-1} \\ K \subset \Omega_k}} \sum_{\substack{\widehat{K} \in \mathcal{M}_k \\ \widehat{K} \subset K}} \sum_{\hat{\sigma} \in \mathcal{F}_{\widehat{K}}} \frac{1}{d_{\widehat{K},\hat{\sigma}}} \int_{\hat{\sigma}} [\varphi_k]_{\widehat{K},\hat{\sigma},\xi_k}(\boldsymbol{y})^2 d\boldsymbol{y}.$$

Let $K$, $\widehat{K}$ and $\hat{\sigma}$ be as in the above sum, either $\hat{\sigma} \subset K$ and so $[\varphi_k]_{\widehat{K},\hat{\sigma},\xi_k} = 0$ or there exists $\sigma \in \mathcal{F}_K$ such that $\hat{\sigma} \subseteq \sigma$. In that latter case, if $(\hat{\sigma},\sigma) \in \widehat{\mathcal{F}}_{k,b} \times \widehat{\mathcal{F}}_{k-1,b}$ or $\hat{\sigma} \in \widehat{\mathcal{F}}_{k,i}$ then $[\varphi_k]_{\widehat{K},\hat{\sigma},\xi_k} = [\varphi_{k-1}]_{K,\sigma,\xi_{k-1}}$. If instead $(\hat{\sigma},\sigma) \in \widehat{\mathcal{F}}_{k,b} \times \widehat{\mathcal{F}}_{k-1,i}$ then $[\varphi_k]_{\widehat{K},\hat{\sigma},\xi_k} = \omega_{K,\sigma}^{-1}[\varphi_{k-1}]_{K,\sigma,\xi_{k-1}}$. See Section 3.1 for an illustration of the above cases. Since $\omega_{K,\sigma}^{-1} \geq 1$, we obtain in all cases

$$|[\varphi_k]_{\widehat{K},\hat{\sigma},\xi_k}| \leq \omega_{K,\sigma}^{-1}|[\varphi_{k-1}]_{K,\sigma,\xi_{k-1}}| \leq C_{\omega,k}|[\varphi_{k-1}]_{K,\sigma,\xi_{k-1}}|.$$

Furthermore, by Assumption 2.18 we have $d_{K,\sigma} \leq C_{\mathrm{r}} d_{\widehat{K},\hat{\sigma}}$. These considerations together give

$$\sum_{\substack{\widehat{K} \in \mathcal{M}_k \\ \widehat{K} \subset K}} \sum_{\hat{\sigma} \in \mathcal{F}_{\widehat{K}}} \frac{1}{d_{\widehat{K},\hat{\sigma}}} \int_{\hat{\sigma}} [\varphi_k]_{\widehat{K},\hat{\sigma},\xi_k}(\boldsymbol{y})^2 d\boldsymbol{y} \leq C_{\omega,k}^2 C_{\mathrm{r}} \sum_{\sigma \in \mathcal{F}_K} \frac{1}{d_{K,\sigma}} \int_\sigma [\varphi_{k-1}]_{K,\sigma,\xi_{k-1}}(\boldsymbol{y})^2 d\boldsymbol{y},$$

hence $|\varphi_k|^2_{\widehat{J}(k),\xi_k} \leq C^2_{\omega,k} C_r |\varphi_{k-1}|^2_{\widehat{J}(k-1),\xi_{k-1}}$ and Eq. (3.2a) is proved. In [11, section 6.1] it is shown that $C_\psi \geq 1$, hence $C_{\omega,k} C_r^{1/2} < C_i = \sqrt{2} C_\psi (1 + C^2_{\omega,k} C_r)^{1/2}$ and Eq. (3.2b) follows. □

The next lemma has been proved for the DGGD scheme in [11] and is valid for the local SWDGGD method as well.

**Lemma 3.2.** *Let $\widehat{\mathcal{D}}_k$ be a local SWDGGD scheme, then there exists $C_{eq} > 0$ depending only on $\alpha, \ell$ and $d$ such that*

$$|\varphi_k|_{\widehat{J}(k),0} \leq C_{eq} \|\nabla_{\widehat{\mathcal{D}}_k} \varphi_k\|_{L^2(\Omega_k)^d} \qquad \text{for all } \varphi_k \in Y_{\mathcal{D}_k}.$$

*Proof.* Follows the lines of [11, Lemma 3.8]. □

The next lemma shows that the error of the local solution depends as usual on the regularity of the solution and data but also on the error committed on the artificial boundary condition. The proof is inspired from the one given in [9] leading to Eq. (2.4) and uses Lemma 2.22.

**Lemma 3.3.** *Let $u \in H^1_0(\Omega)$ be the exact solution to Eq. (2.2), $\kappa_k \in Z_{\mathcal{D}_k}$ and $\hat{\vartheta}_k \in Y_{\mathcal{D}_k}$ be solution of Eq. (3.1c). Then*

$$\begin{aligned}
&\|\nabla u - \nabla_{\widehat{\mathcal{D}}_k,\kappa_k} \hat{\vartheta}_k\|_{L^2(\Omega_k)^d} + |\hat{\vartheta}_k|_{\widehat{J}(k),\kappa_k} \\
&\leq \frac{1 + C_{eq}}{\underline{\lambda}} W_{\widehat{\mathcal{D}}_k}(A \nabla u + \boldsymbol{F}) + C_A \min_{\xi_k \in Z_{\mathcal{D}_k}} (S_{\widehat{\mathcal{D}}_k,J,\xi_k}(u) + C_\partial |\kappa_k - \xi_k|_{\partial \Omega_k^-})
\end{aligned} \qquad (3.3)$$

*with $C_A := C_{eq}(1 + \kappa(A))$ and $C_\partial$ from Lemma 2.22.*

**Remark 3.4.** *Observe that Lemma 3.3 is valid for any $\kappa_k \in Z_{\mathcal{D}_k}$ and not only $\kappa_k$ given by scheme Eq. (3.1).*

*Proof.* Since $\widehat{\mathcal{D}}_k$ is a SWDGGD scheme, by Definition 2.7 for any $\boldsymbol{v} \in H_{\text{div}}(\Omega_k)$ and $\psi_k \in Y_{\mathcal{D}_k}$ we have

$$\left| \int_{\Omega_k} (\nabla_{\widehat{\mathcal{D}}_k} \psi_k \cdot \boldsymbol{v} + \Pi_{\widehat{\mathcal{D}}_k} \psi_k \nabla \cdot \boldsymbol{v}) \, d\boldsymbol{x} \right| \leq \|\nabla_{\widehat{\mathcal{D}}_k} \psi_k\|_{L^2(\Omega_k)^d} W_{\widehat{\mathcal{D}}_k}(\boldsymbol{v}).$$

As $-\nabla \cdot (A \nabla u + \boldsymbol{F}) = f_0 \in L^2(\Omega_k)$ we can take $\boldsymbol{v} = A \nabla u + \boldsymbol{F}$ and obtain

$$\left| \int_{\Omega_k} (\nabla_{\widehat{\mathcal{D}}_k} \psi_k \cdot (A \nabla u + \boldsymbol{F}) - \Pi_{\widehat{\mathcal{D}}_k} \psi_k f_0) \, d\boldsymbol{x} \right| \leq \|\nabla_{\widehat{\mathcal{D}}_k} \psi_k\|_{L^2(\Omega_k)^d} W_{\widehat{\mathcal{D}}_k}(A \nabla u + \boldsymbol{F}).$$

Using Eq. (3.1c) we get

$$\left| \int_{\Omega_k} A(\nabla u - \nabla_{\widehat{\mathcal{D}}_k,\kappa_k} \hat{\vartheta}_k) \cdot \nabla_{\widehat{\mathcal{D}}_k} \psi_k \, d\boldsymbol{x} \right| \leq \|\nabla_{\widehat{\mathcal{D}}_k} \psi_k\|_{L^2(\Omega_k)^d} W_{\widehat{\mathcal{D}}_k}(A \nabla u + \boldsymbol{F}).$$

Let $\varphi_k \in Y_{\mathcal{D}_k}$, we have

$$\int_{\Omega_k} A(\nabla_{\widehat{\mathcal{D}}_k,\kappa_k} \varphi_k - \nabla_{\widehat{\mathcal{D}}_k,\kappa_k} \hat{\vartheta}_k) \cdot \nabla_{\widehat{\mathcal{D}}_k} \psi_k \, d\boldsymbol{x}$$

$$\leq \|\nabla_{\widehat{\mathcal{D}}_k} \psi_k\|_{L^2(\Omega_k)^d} W_{\widehat{\mathcal{D}}_k}(A \nabla u + \boldsymbol{F}) + \int_{\Omega_k} A(\nabla_{\widehat{\mathcal{D}}_k,\kappa_k} \varphi_k - \nabla u) \cdot \nabla_{\widehat{\mathcal{D}}_k} \psi_k \, d\boldsymbol{x}$$

$$\leq \|\nabla_{\widehat{\mathcal{D}}_k} \psi_k\|_{L^2(\Omega_k)^d} (W_{\widehat{\mathcal{D}}_k}(A \nabla u + \boldsymbol{F}) + \overline{\lambda} \|\nabla_{\widehat{\mathcal{D}}_k,\kappa_k} \varphi_k - \nabla u\|_{L^2(\Omega_k)^d}).$$

We choose $\psi_k = \varphi_k - \hat{\vartheta}_k$, since $\nabla_{\widehat{\mathcal{D}}_k, \kappa_k} \varphi_k - \nabla_{\widehat{\mathcal{D}}_k, \kappa_k} \hat{\vartheta}_k = \nabla_{\widehat{\mathcal{D}}_k, 0}(\varphi_k - \hat{\vartheta}_k) = \nabla_{\widehat{\mathcal{D}}_k}(\varphi_k - \hat{\vartheta}_k)$ we get

$$\underline{\lambda} \|\nabla_{\widehat{\mathcal{D}}_k}(\varphi_k - \hat{\vartheta}_k)\|^2_{L^2(\Omega_k)^d} \leq \int_{\Omega_k} A(\nabla_{\widehat{\mathcal{D}}_k, \kappa_k} \varphi_k - \nabla_{\widehat{\mathcal{D}}_k, \kappa_k} \hat{\vartheta}_k) \cdot \nabla_{\widehat{\mathcal{D}}_k}(\varphi_k - \hat{\vartheta}_k) \, d\boldsymbol{x}$$

and hence

$$\|\nabla_{\widehat{\mathcal{D}}_k}(\varphi_k - \hat{\vartheta}_k)\|_{L^2(\Omega_k)^d} \leq \frac{1}{\underline{\lambda}} W_{\widehat{\mathcal{D}}_k}(A\nabla u + \boldsymbol{F}) + \kappa(A)\|\nabla_{\widehat{\mathcal{D}}_k, \kappa_k} \varphi_k - \nabla u\|_{L^2(\Omega_k)^d}. \tag{3.4}$$

This gives

$$\begin{aligned}\|\nabla u - \nabla_{\widehat{\mathcal{D}}_k, \kappa_k} \hat{\vartheta}_k\|_{L^2(\Omega_k)^d} \leq &\|\nabla u - \nabla_{\widehat{\mathcal{D}}_k, \kappa_k} \varphi_k\|_{L^2(\Omega_k)^d} + \|\nabla_{\widehat{\mathcal{D}}_k}(\varphi_k - \hat{\vartheta}_k)\|_{L^2(\Omega_k)^d} \\ \leq &\frac{1}{\underline{\lambda}} W_{\widehat{\mathcal{D}}_k}(A\nabla u + \boldsymbol{F}) + (1 + \kappa(A)) \|\nabla u - \nabla_{\widehat{\mathcal{D}}_k, \kappa_k} \varphi_k\|_{L^2(\Omega_k)^d}. \end{aligned} \tag{3.5}$$

Using $|\hat{\vartheta}_k|_{\widehat{J}(k), \kappa_k} \leq |\varphi_k|_{\widehat{J}(k), \kappa_k} + |\hat{\vartheta}_k - \varphi_k|_{\widehat{J}(k), 0}$, Lemma 3.2 and Eq. (3.4) yields

$$\begin{aligned}|\hat{\vartheta}_k|_{\widehat{J}(k), \kappa_k} \leq &|\varphi_k|_{\widehat{J}(k), \kappa_k} + \frac{C_{\mathrm{eq}}}{\underline{\lambda}} W_{\widehat{\mathcal{D}}_k}(A\nabla u + \boldsymbol{F}) \\ &+ C_{\mathrm{eq}} \kappa(A) \|\nabla_{\widehat{\mathcal{D}}_k, \kappa_k} \varphi_k - \nabla u\|_{L^2(\Omega_k)^d}. \end{aligned} \tag{3.6}$$

Summing Eqs. (3.5) and (3.6) and taking the infimum over $\varphi_k \in Y_{\mathcal{D}_k}$ we get

$$\begin{aligned}\|\nabla u - \nabla_{\widehat{\mathcal{D}}_k, \kappa_k} \hat{\vartheta}_k\|_{L^2(\Omega_k)^d} + |\hat{\vartheta}_k|_{\widehat{J}(k), \kappa_k} \leq &\frac{1 + C_{\mathrm{eq}}}{\underline{\lambda}} W_{\widehat{\mathcal{D}}_k}(A\nabla u + \boldsymbol{F}) \\ &+ C_{\mathrm{eq}}(1 + \kappa(A)) S_{\widehat{\mathcal{D}}_k, J, \kappa_k}(u). \end{aligned}$$

We conclude using Lemma 2.22 and taking the inf over $\xi_k$. $\square$

**Lemma 3.5.** *Let $((\kappa_k, \hat{\vartheta}_k))_{k=1}^M$ be the sequence defined by the local elliptic scheme Eqs. (3.1a) to (3.1c). Then for $k \geq 2$*

$$\min_{\xi_k \in Z_{\mathcal{D}_k}} (S_{\widehat{\mathcal{D}}_k, J, \xi_k}(u) + C_\partial |\kappa_k - \xi_k|_{\partial \Omega_k^-})$$
$$\leq 2C_{\mathrm{i}} \left( \|\nabla_{\widehat{\mathcal{D}}_k, \kappa_{k-1}} \hat{\vartheta}_{k-1} - \nabla u\|_{L^2(\Omega_k)^d} + |\hat{\vartheta}_{k-1}|_{\widehat{J}(k-1), \kappa_{k-1}} \right),$$

*where $C_{\mathrm{i}}$ is defined in Lemma 3.1.*

*Proof.* Taking $\xi_k = \kappa_k$ we have

$$\min_{\xi_k \in Z_{\mathcal{D}_k}} (S_{\widehat{\mathcal{D}}_k, J, \xi_k}(u) + C_\partial |\kappa_k - \xi_k|_{\partial \Omega_k^-}) \leq S_{\widehat{\mathcal{D}}_k, J, \kappa_k}(u).$$

Since $\kappa_k = (\kappa_{k-1} + \hat{\vartheta}_{k-1}) \chi_{\Omega \setminus \Omega_k}$ by Lemma 3.1 there exists $\hat{\varphi}_k \in Y_{\mathcal{D}_k}$ satisfying

$$\|\nabla_{\widehat{\mathcal{D}}_k, \kappa_k} \hat{\varphi}_k - \nabla_{\widehat{\mathcal{D}}_{k-1}, \kappa_{k-1}} \hat{\vartheta}_{k-1}\|_{L^2(\Omega_k)^d} \leq C_{\mathrm{i}} |\hat{\vartheta}_{k-1}|_{\widehat{J}(k-1), \kappa_{k-1}},$$
$$|\hat{\varphi}_k|_{\widehat{J}(k), \kappa_k} \leq C_{\mathrm{i}} |\hat{\vartheta}_{k-1}|_{\widehat{J}(k-1), \kappa_{k-1}}$$

and so

$$S_{\widehat{\mathcal{D}}_k, J, \kappa_k}(u) = \inf_{\varphi \in Y_{\mathcal{D}_k}} (\|\nabla_{\widehat{\mathcal{D}}_k, \kappa_k} \varphi - \nabla u\|_{L^2(\Omega_k)^d} + |\varphi|_{\widehat{J}(k), \kappa_k})$$
$$\leq \|\nabla_{\widehat{\mathcal{D}}_k, \kappa_k} \hat{\varphi}_k - \nabla u\|_{L^2(\Omega_k)^d} + |\hat{\varphi}_k|_{\widehat{J}(k), \kappa_k}$$
$$\leq \|\nabla_{\widehat{\mathcal{D}}_k, \kappa_{k-1}} \hat{\vartheta}_{k-1} - \nabla u\|_{L^2(\Omega_k)^d} + 2C_i |\hat{\vartheta}_{k-1}|_{\widehat{J}(k-1), \kappa_{k-1}}.$$

□

Let $\mathcal{H} \subset \mathbb{R}_+$ be a countable set with 0 as only accumulation point. For each $h \in \mathcal{H}$ we consider a polytopal mesh sequence $(\mathfrak{T}_{h,k})_{k=1}^M = ((\mathcal{M}_{h,k}, \mathcal{F}_{h,k}, \mathcal{P}_{h,k}))_{k=1}^M$ satisfying Assumption 2.18 with $h = \max_{k=1,...,M} h_{\mathcal{M}_{h,k}}$, where $h_{\mathcal{M}_{h,k}} = \max\{h_K : K \in \mathcal{M}_{h,k}\}$. Let $\mathcal{D}_{h,k}$ and $\widehat{\mathcal{D}}_{h,k}$ be the global and local SWDGGD schemes given by those meshes $\mathfrak{T}_{h,k}$. In the following the index $h$ in $\mathcal{D}_{h,k}$ and $\widehat{\mathcal{D}}_{h,k}$ is left out of notation for the sake of simplicity.

**Theorem 3.6.** *Let $\mathcal{D}_k$ and $\widehat{\mathcal{D}}_k$ be global and local SWDGGD. Let $((\kappa_k, \hat{\vartheta}_k))_{k=1}^M$ be the sequence defined by the local elliptic scheme Eqs. (3.1a) to (3.1c) and $u \in H_0^1(\Omega)$ the exact solution to Eq. (2.2). Then for $k = 1, ..., M$*

$$\lim_{h \to 0} \|\nabla u - \nabla_{\widehat{\mathcal{D}}_k, \kappa_k} \hat{\vartheta}_k\|_{L^2(\Omega_k)^d} + |\hat{\vartheta}_k|_{\widehat{J}(k), \kappa_k} = 0. \tag{3.7a}$$

*If moreover $u \in H_0^1(\Omega) \cap H^2(\Omega)$, the coefficients of $A$ are Lipschitz continuous and $\boldsymbol{F} \in H^1(\Omega)^d$ there exists $C_1, C_2, C_3$ depending on $\alpha, \ell, d, \rho, C_r, |\Omega|, A, \boldsymbol{F}$ and $u$ such that*

$$\|\nabla u - \nabla_{\widehat{\mathcal{D}}_k, \kappa_k} \hat{\vartheta}_k\|_{L^2(\Omega_k)^d} + |\hat{\vartheta}_k|_{\widehat{J}(k), \kappa_k} \leq C_1 h, \tag{3.7b}$$

$$\|\nabla u - \nabla_{\widehat{\mathcal{D}}_k, \kappa_k} \hat{\vartheta}_k\|_{L^2(\Omega_k)^d} + |\hat{\vartheta}_k|_{\widehat{J}(k), \kappa_k} \leq C_2 h_{\widehat{\mathcal{M}}_k} + C_3 |\kappa_k - \xi_k|_{\partial \Omega_k^-}, \tag{3.7c}$$

*where $\xi_k = \operatorname{argmin}_{\xi \in Z_{\mathcal{D}_k}} S_{\widehat{\mathcal{D}}_k, J, \xi}(u)$.*

**Remark 3.7.** *The above theorem gives three important results. The first one Eq. (3.7a) asserts that the numerical solution given by the local scheme Eqs. (3.1a) to (3.1c) converges to the exact solution even under weak regularity of the solution and data. Assuming more regularity we recover in Eq. (3.7b) the usual convergence rate. In Eq. (3.7c) we establish that the error on the local domain depends on the local mesh size and the error committed on the artificial boundary.*

*Proof of Theorem 3.6.* Let

$$E_{\widehat{\mathcal{D}}_k} := \|\nabla u - \nabla_{\widehat{\mathcal{D}}_k, \kappa_k} \hat{\vartheta}_k\|_{L^2(\Omega_k)^d} + |\hat{\vartheta}_k|_{\widehat{J}(k), \kappa_k}.$$

Since $\kappa_1 = 0$ and $Z_{\mathcal{D}_1} = \{0\}$ by Lemma 3.3 we have

$$E_{\widehat{\mathcal{D}}_1} \leq \frac{1 + C_{\text{eq}}}{\underline{\lambda}} W_{\widehat{\mathcal{D}}_1}(A \nabla u + \boldsymbol{F}) + C_A S_{\widehat{\mathcal{D}}_1, J, 0}(u)$$

and by Lemmas 3.3 and 3.5 we have, for $k \geq 2$,

$$E_{\widehat{\mathcal{D}}_k} \leq \frac{1 + C_{\text{eq}}}{\underline{\lambda}} W_{\widehat{\mathcal{D}}_k}(A \nabla u + \boldsymbol{F}) + 2C_i C_A E_{\widehat{\mathcal{D}}_{k-1}}.$$

Let $\alpha = 2C_iC_A$, since $S_{\widehat{\mathcal{D}}_1,J,0}(u) \leq S_{\mathcal{D}_1,J}(u)$ it holds

$$E_{\widehat{\mathcal{D}}_k} \leq \alpha^{k-1} E_{\widehat{\mathcal{D}}_1} + \frac{1+C_{\text{eq}}}{\underline{\lambda}} \sum_{j=2}^{k} \alpha^{k-j} W_{\widehat{\mathcal{D}}_j}(A\nabla u + \boldsymbol{F})$$
$$\leq C_A \alpha^{k-1} S_{\mathcal{D}_1,J}(u) + \frac{1+C_{\text{eq}}}{\underline{\lambda}} \sum_{j=1}^{k} \alpha^{k-j} W_{\widehat{\mathcal{D}}_j}(A\nabla u + \boldsymbol{F}). \tag{3.8}$$

We have thus proved Eq. (3.7a) thanks to Eq. (3.8), Lemma 2.16 and the limit conformity of $\widehat{\mathcal{D}}_j$ for $j = 1, ..., k$ (we recall that $\widehat{\mathcal{D}}_j$ is a SWDGGD and hence a sequence of $\widehat{\mathcal{D}}_j$ is limit conforming). Under the additional assumptions on the data, from Lemmas 2.13 and 2.14 for $\widehat{\mathcal{D}}_k$ we have

$$S_{\mathcal{D}_1,J}(u) \leq C_S h_{\mathcal{M}_1} \|u\|_{H^2(\Omega)},$$
$$W_{\widehat{\mathcal{D}}_k}(A\nabla u + \boldsymbol{F}) \leq C_W h_{\widehat{\mathcal{M}}_k} \|A\nabla u + \boldsymbol{F}\|_{H^1(\Omega_k)^d} \tag{3.9}$$

and so

$$E_{\widehat{\mathcal{D}}_k} \leq C_A \alpha^{k-1} C_S h_{\mathcal{M}_1} \|u\|_{H^2(\Omega)} + \frac{1+C_{\text{eq}}}{\underline{\lambda}} \sum_{j=1}^{k} \alpha^{k-1} C_W h_{\widehat{\mathcal{M}}_j} \|A\nabla u + \boldsymbol{F}\|_{H^1(\Omega_j)^d},$$

which implies Eq. (3.7b) with

$$C_1 := C_A \alpha^{k-1} C_S \|u\|_{H^2(\Omega)} + \frac{1+C_{\text{eq}}}{\underline{\lambda}} \sum_{j=1}^{k} \alpha^{k-1} C_W \|A\nabla u + \boldsymbol{F}\|_{H^1(\Omega_j)^d}.$$

Let $\xi_k = \operatorname{argmin}_{\xi \in Z_{\mathcal{D}_k}} S_{\widehat{\mathcal{D}}_k,J,\xi}(u)$, it holds

$$\min_{\xi \in Z_{\mathcal{D}_k}} (S_{\widehat{\mathcal{D}}_k,J,\xi}(u) + C_\partial |\kappa_k - \xi|_{\partial\Omega_k^-}) \leq S_{\widehat{\mathcal{D}}_k,J,\xi_k}(u) + C_\partial |\kappa_k - \xi_k|_{\partial\Omega_k^-},$$

using Lemma 2.21 we get $S_{\widehat{\mathcal{D}}_k,J,\xi_k}(u) \leq C_S h_{\widehat{\mathcal{M}}_k} \|u\|_{H^2(\Omega_k)}$ and again from Lemma 3.3 and Eq. (3.9) we obtain the bound Eq. (3.7c) with

$$C_2 := C_A C_S \|u\|_{H^2(\Omega_k)} + \frac{C_W}{\underline{\lambda}} \|A\nabla u + \boldsymbol{F}\|_{H^1(\Omega_k)^d}, \tag{3.10}$$

where $C_3 = C_A C_\partial$. □

## 3.2 Improved local estimate

Under stronger conditions and using the pointwise error estimates proved in [6] we can further improve the local estimate Eq. (3.7c) for $k = 2$.

Let $\boldsymbol{z} \in \Omega$, the weight function $\sigma_{\boldsymbol{z},h}(\boldsymbol{x}) = h/(h + |\boldsymbol{x} - \boldsymbol{z}|)$ and $\|\cdot\|_{W^{2,\infty}_{\boldsymbol{z},h}(\Omega)}$ a weighted Sobolev norm defined by

$$\|v\|_{W^{2,\infty}_{\boldsymbol{z},h}(\Omega)} = \max_{i=0,1,2} |v|_{W^{i,\infty}_{\boldsymbol{z},h}}, \qquad |v|_{W^{i,\infty}_{\boldsymbol{z},h}} = \max_{|\alpha|=i} \|\sigma_{\boldsymbol{z},h} \frac{\partial^\alpha v}{\partial \boldsymbol{x}^\alpha}\|_{L^\infty(\Omega)}.$$

We will use the following lemma, which is a version of [6, Corollary 5.5].

**Lemma 3.8.** Let $A = aI_d$ with $I_d \in \mathbb{R}^{d \times d}$ the identity matrix and $a > 0$. Let $u \in W_0^{1,\infty}(\Omega) \cap W^{2,\infty}(\Omega)$ be solution of Eq. (2.1) with $f \in L^2(\Omega)$, $\vartheta_1 \in X_{\mathcal{D}_1}$ solution of Eq. (3.1c). Then there is a constant $C_\infty >$ such that for any $\boldsymbol{z} \in \overline{\Omega}$

$$|u(\boldsymbol{z}) - \Pi_{\mathcal{D}_1}\vartheta_1(\boldsymbol{z})| \leq C_\infty h^2 \log(h^{-1}) \|u\|_{W_{\boldsymbol{z},h}^{2,\infty}(\Omega)}. \tag{3.11}$$

Applying Lemma 3.8 to Eq. (3.7c) we obtain a better bound on the local error for $k = 2$, as explained in the following theorem.

**Theorem 3.9.** Let $u \in W_0^{1,\infty}(\Omega) \cap W^{2,\infty}(\Omega)$ be solution of Eq. (2.1) with $A = aI_d$ and $f \in L^2(\Omega)$ as in Lemma 3.8. Let $\mathcal{D}_k$ and $\widehat{\mathcal{D}}_k$ be global and local SWDGGD schemes and $((\kappa_k, \hat{\vartheta}_k))_{k=1}^2$ the sequence defined by the local elliptic scheme Eqs. (3.1a) to (3.1c). Under the assumption that $h \leq C_h \min_{K \in \mathcal{M}_1} h_K$ for $C_h > 0$ it exists $C_4$ independent of $u$ and $h$ such that

$$\begin{aligned}
\|\nabla u - \nabla_{\widehat{\mathcal{D}}_2, \kappa_2} \hat{\vartheta}_2\|_{L^2(\Omega_2)^d} &+ |\hat{\vartheta}_2|_{\widehat{J}(2), \kappa_2} \\
&\leq C_2 h_{\widehat{\mathcal{M}}_2} + C_4 (h_{\widehat{\mathcal{M}}_2} |u|_{H^2(D_2)} + h^{3/2} \log(h^{-1}) \sup_{\boldsymbol{y} \in \partial \Omega_2 \setminus \partial \Omega} \|u\|_{W_{\boldsymbol{y},h}^{2,\infty}(\Omega)}),
\end{aligned} \tag{3.12}$$

where $D_2$ is a neighborhood of $\Omega_2$ specified below.

**Remark 3.10.** Equation (3.12) bounds the error in the local domain $\Omega_2$ and has three terms in the right. From Eq. (3.10) we see that the first term depends on $u$ and $\boldsymbol{F}$ in $\Omega_2$. The second term depends on $u$ in a small neighborhood of $\Omega_2$. The last term depends on the regularity of $u$ in the whole domain, but it is of higher order and is measured in a weighted norm which weight is $\mathcal{O}(1)$ close to the artificial boundary and $\mathcal{O}(h)$ far from it. Hence the error in $\Omega_2$ depends mostly on the regularity of $u$ and $\boldsymbol{F}$ inside or very close to $\Omega_2$.

*Proof.* First we observe that Eq. (3.7c) for $k = 2$ is valid with $\xi_2 \in Z_{\mathcal{D}_2}$ such that $\Pi_{\mathcal{D}_2}\xi_2$ is the orthogonal projection of $u$ onto $\Pi_{\mathcal{D}_2} Z_{\mathcal{D}_2}$, indeed even for this choice of $\xi_2$ we still have $S_{\widehat{\mathcal{D}}_2, J, \xi_2}(u) \leq C_S h_{\widehat{\mathcal{M}}_2} \|u\|_{H^2(\Omega_2)}$. Let $K \in \widehat{\mathcal{M}}_2$, $T \in \mathcal{M}_2 \setminus \widehat{\mathcal{M}}_2$ and $\sigma \in \mathcal{F}_K \cap \mathcal{F}_T$. From Assumption 2.18b) we have $K, T \in \mathcal{M}_1$ and Assumption 2.18d) implies $h_K \leq \rho h_T$. There exists $C_\Pi$ ([7, Lemma 1.59]) independent of $u$, $T$ and $h_K$ such that

$$\int_\sigma |\Pi_{\overline{T}} \xi_2 - u|(\boldsymbol{y})^2 \, d\boldsymbol{y} \leq C_\Pi h_K^3 |u|_{H^2(T)}^2.$$

Using Assumption 2.18d) we obtain $1/d_{K,\sigma} \leq \rho/h_K$, hence

$$\sum_{\substack{K \in \widehat{\mathcal{M}}_2 \\ T \in \mathcal{M}_2 \setminus \widehat{\mathcal{M}}_2}} \sum_{\sigma \in \mathcal{F}_K \cap \mathcal{F}_T} \frac{1}{d_{K,\sigma}} \int_\sigma |\Pi_{\overline{T}} \xi_2 - u|(\boldsymbol{y})^2 \, d\boldsymbol{y} \leq C_\Pi \rho h_{\widehat{\mathcal{M}}_2}^2 |u|_{H^2(D_2)}^2,$$

with $D_2 = \cup_{\{T \in \mathcal{M}_2 \setminus \widehat{\mathcal{M}}_2 : \partial T \cap \partial K \neq \emptyset, K \in \widehat{\mathcal{M}}_2\}} T$. From Lemma 3.8 we have

$$\int_\sigma |u - \Pi_{\overline{T}} \kappa_2|(\boldsymbol{y})^2 \, d\boldsymbol{y} = \int_\sigma |u - \Pi_{\overline{T}} \vartheta_1|(\boldsymbol{y})^2 \, d\boldsymbol{y} \leq |\sigma| C_\infty^2 h^4 \log(h^{-1})^2 \sup_{\boldsymbol{y} \in \sigma} \|u\|_{W_{\boldsymbol{y},h}^{2,\infty}(\Omega)}^2.$$

Since $h \leq C_h \min_{K \in \mathcal{M}_1} h_K$ it follows that $1/d_{K,\sigma} \leq C_h \rho/h$ and thus

$$\sum_{\substack{K \in \widehat{\mathcal{M}}_2 \\ T \in \mathcal{M}_2 \setminus \widehat{\mathcal{M}}_2}} \sum_{\sigma \in \mathcal{F}_K \cap \mathcal{F}_T} \frac{1}{d_{K,\sigma}} \int_\sigma |u - \Pi_{\overline{T}} \kappa_2|(\boldsymbol{y})^2 \, d\boldsymbol{y}$$

$$\leq |\partial \Omega_2 \setminus \partial \Omega| C_\infty^2 C_h \rho h^3 \log(h^{-1})^2 \sup_{\boldsymbol{y} \in \partial \Omega_2 \setminus \partial \Omega} \|u\|_{W_{\boldsymbol{y},h}^{2,\infty}(\Omega)}^2.$$

Applying a triangle inequality on $|\kappa_2 - \xi_2|_{\partial \Omega_2^-}$ in Eq. (3.7c) we get Eq. (3.12). □

## 3.3 A priori error analysis for the global solution

We next study the error on the whole domain $\Omega$ of the numerical solution $\vartheta_k \in X_{\mathcal{D}_k}$ defined by our local scheme Eq. (3.1). The next Lemma 3.11 is the main ingredient for the global error bound.

**Lemma 3.11.** *Let $u \in H_0^1(\Omega)$ be solution to Eq. (2.2) and $(\vartheta_k)_{k=1}^M$ be the sequence defined by scheme Eqs. (3.1a) to (3.1c). Then we have*

$$\|\nabla u - \nabla_{\mathcal{D}_k}\vartheta_k\|_{L^2(\Omega)^d} + |\vartheta_k|_{J(k)} \leq C_5(\|\nabla u - \nabla_{\mathcal{D}_{k-1}}\vartheta_{k-1}\|_{L^2(\Omega)^d} + |\vartheta_{k-1}|_{J(k-1)})$$
$$+ C_5(\|\nabla u - \nabla_{\widehat{\mathcal{D}}_k, \kappa_k}\hat{\vartheta}_k\|_{L^2(\Omega_k)^d} + |\hat{\vartheta}_k|_{\widehat{J}(k),\kappa_k}).$$

*where $C_5 = \sqrt{2}(1 + C_\psi)(1 + \sqrt{2}C_{\omega,k})$.*

*Proof.* We have

$$\|\nabla u - \nabla_{\mathcal{D}_k}\vartheta_k\|_{L^2(\Omega)^d}^2$$
$$= \sum_{K \in \mathcal{M}_k} \sum_{\sigma \in \mathcal{F}_K} \int_{D_{K,\sigma}} |\nabla u(\boldsymbol{x}) - \nabla_{\overline{K}}\vartheta_k(\boldsymbol{x}) - \psi(s)\frac{[\vartheta_k]_{K,\sigma}(\boldsymbol{y})}{d_{K,\sigma}}\boldsymbol{n}_{K,\sigma}|^2 d\boldsymbol{x}$$
$$= \sum_{T \in \mathcal{M}_k \setminus \widehat{\mathcal{M}}_k} \sum_{\sigma \in \mathcal{F}_T} \int_{D_{T,\sigma}} |\nabla u(\boldsymbol{x}) - \nabla_{\overline{T}}\vartheta_k(\boldsymbol{x}) - \psi(s)\frac{[\vartheta_k]_{T,\sigma}(\boldsymbol{y})}{d_{T,\sigma}}\boldsymbol{n}_{T,\sigma}|^2 d\boldsymbol{x}$$
$$+ \sum_{K \in \widehat{\mathcal{M}}_k} \sum_{\sigma \in \mathcal{F}_K} \int_{D_{K,\sigma}} |\nabla u(\boldsymbol{x}) - \nabla_{\overline{K}}\vartheta_k(\boldsymbol{x}) - \psi(s)\frac{[\vartheta_k]_{K,\sigma}(\boldsymbol{y})}{d_{K,\sigma}}\boldsymbol{n}_{K,\sigma}|^2 d\boldsymbol{x}$$
$$= I + II.$$

For the first term $I$, we have the following considerations. Let $T \in \mathcal{M}_k \setminus \widehat{\mathcal{M}}_k$, then $T \in \mathcal{M}_{k-1}$ and $\nabla_{\overline{T}}\vartheta_k = \nabla_{\overline{T}}\vartheta_{k-1}$. Let $\sigma \in \mathcal{F}_T$, if $\sigma \notin \widehat{\mathcal{F}}_{k,b}$ then $[\vartheta_k]_{T,\sigma} = [\vartheta_{k-1}]_{T,\sigma}$. If $\sigma \in \widehat{\mathcal{F}}_{k,b}$ then $\sigma = \partial K \cap \partial T$ with $K \in \widehat{\mathcal{M}}_k$ and by Assumption 2.18b) $K \in \widehat{\mathcal{M}}_{k-1}$. Using Eqs. (3.1a) and (3.1b) we have

$$[\vartheta_k]_{T,\sigma} - [\vartheta_{k-1}]_{T,\sigma} = \omega_{T,\sigma}(\Pi_{\overline{K}}\hat{\vartheta}_k - \Pi_{\overline{T}}\vartheta_{k-1}) - \omega_{T,\sigma}(\Pi_{\overline{K}}\vartheta_{k-1} - \Pi_{\overline{T}}\vartheta_{k-1})$$
$$= \omega_{T,\sigma}(\Pi_{\overline{K}}\hat{\vartheta}_k - \Pi_{\overline{K}}\vartheta_{k-1}).$$

Next, adding and removing $[\vartheta_{k-1}]_{T,\sigma}$ from $[\vartheta_k]_{T,\sigma}$ we get

$$I \leq 2 \sum_{T \in \mathcal{M}_k \setminus \widehat{\mathcal{M}}_k} \sum_{\sigma \in \mathcal{F}_T} \int_{D_{T,\sigma}} |\nabla u(\boldsymbol{x}) - \nabla_{\overline{T}}\vartheta_{k-1}(\boldsymbol{x}) - \psi(s)\frac{[\vartheta_{k-1}]_{T,\sigma}(\boldsymbol{y})}{d_{T,\sigma}}\boldsymbol{n}_{T,\sigma}|^2 d\boldsymbol{x}$$
$$+ 2 \sum_{\substack{K \in \widehat{\mathcal{M}}_k \\ T \in \mathcal{M}_k \setminus \widehat{\mathcal{M}}_k}} \sum_{\sigma \in \mathcal{F}_T \cap \mathcal{F}_K} \int_{D_{T,\sigma}} |\psi(s)\omega_{T,\sigma}\frac{(\Pi_{\overline{K}}\hat{\vartheta}_k - \Pi_{\overline{K}}\vartheta_{k-1})(\boldsymbol{y})}{d_{T,\sigma}}|^2 d\boldsymbol{x}.$$

Since $\omega_{T,\sigma} \leq 1$, using a change of variables we have

$$
\begin{aligned}
I \leq & 2\|\nabla u - \nabla_{\mathcal{D}_{k-1}}\vartheta_{k-1}\|^2_{L^2(\Omega\setminus\Omega_k)} \\
& + 2C_\psi^2 \sum_{\substack{K\in\widehat{\mathcal{M}}_k \\ T\in\mathcal{M}_k\setminus\widehat{\mathcal{M}}_k}} \sum_{\sigma\in\mathcal{F}_T\cap\mathcal{F}_K} \frac{1}{d_{T,\sigma}} \int_\sigma (\Pi_{\overline{K}}\hat{\vartheta}_k - \Pi_{\overline{K}}\vartheta_{k-1})(\boldsymbol{y})^2 d\boldsymbol{y} \\
= & 2\|\nabla u - \nabla_{\mathcal{D}_{k-1}}\vartheta_{k-1}\|^2_{L^2(\Omega\setminus\Omega_k)} + 2C_\psi^2|\hat{\vartheta}_k - \vartheta_{k-1}|^2_{\partial\Omega_k^+},
\end{aligned}
$$

where for $\phi_k \in X_{\mathcal{D}_k}$

$$
|\phi_k|^2_{\partial\Omega_k^+} := \sum_{\substack{K\in\widehat{\mathcal{M}}_k \\ T\in\mathcal{M}_k\setminus\widehat{\mathcal{M}}_k}} \sum_{\sigma\in\mathcal{F}_K\cap\mathcal{F}_T} \frac{1}{d_{T,\sigma}} \int_\sigma \Pi_{\overline{K}}\phi_k(\boldsymbol{y})^2 \, d\boldsymbol{y}.
$$

For the second term $II$ we have $[\vartheta_k]_{K,\sigma} = [\hat{\vartheta}_k]_{K,\sigma,\kappa_k}$ if $\sigma \in \widehat{\mathcal{F}}_{k,i}$ or $\sigma \in \widehat{\mathcal{F}}_{k,b} \cap \mathcal{F}_{k,b}$ and $[\vartheta_k]_{K,\sigma} = \omega_{K,\sigma}[\hat{\vartheta}_k]_{K,\sigma,\kappa_k}$ if $\sigma \in \widehat{\mathcal{F}}_{k,b} \setminus \mathcal{F}_{k,b}$. Hence

$$
\begin{aligned}
II \leq & 2\|\nabla u - \nabla_{\widehat{\mathcal{D}}_k,\kappa_k}\hat{\vartheta}_k\|^2_{L^2(\Omega_k)^d} \\
& + 2 \sum_{\substack{K\in\widehat{\mathcal{M}}_k \\ T\in\mathcal{M}_k\setminus\widehat{\mathcal{M}}_k}} \sum_{\sigma\in\mathcal{F}_T\cap\mathcal{F}_K} \int_{D_{K,\sigma}} |\psi(s)\frac{(1-\omega_{K,\sigma})[\hat{\vartheta}_k]_{K,\sigma,\kappa_k}(\boldsymbol{y})}{d_{K,\sigma}}|^2 \, d\boldsymbol{x} \\
\leq & 2\|\nabla u - \nabla_{\widehat{\mathcal{D}}_k,\kappa_k}\hat{\vartheta}_k\|^2_{L^2(\Omega_k)^d} + 2C_\psi^2 \sum_{\substack{K\in\widehat{\mathcal{M}}_k \\ T\in\mathcal{M}_k\setminus\widehat{\mathcal{M}}_k}} \sum_{\sigma\in\mathcal{F}_T\cap\mathcal{F}_K} \frac{1}{d_{K,\sigma}} \int_\sigma [\hat{\vartheta}_k]_{K,\sigma,\kappa_k}(\boldsymbol{y})^2 \, d\boldsymbol{y} \\
\leq & 2\|\nabla u - \nabla_{\widehat{\mathcal{D}}_k,\kappa_k}\hat{\vartheta}_k\|^2_{L^2(\Omega_k)^d} + 2C_\psi^2|\hat{\vartheta}_k|^2_{\widehat{J}(k),\kappa_k}.
\end{aligned}
$$

We then obtain

$$
\begin{aligned}
\|\nabla u - \nabla_{\mathcal{D}_k}\vartheta_k\|^2_{L^2(\Omega)^d} \leq & 2\|\nabla u - \nabla_{\mathcal{D}_{k-1}}\vartheta_{k-1}\|^2_{L^2(\Omega)^d} + 2\|\nabla u - \nabla_{\widehat{\mathcal{D}}_k,\kappa_k}\hat{\vartheta}_k\|^2_{L^2(\Omega_k)^d} \\
& + 2C_\psi^2|\hat{\vartheta}_k - \vartheta_{k-1}|^2_{\partial\Omega_k^+} + 2C_\psi^2|\hat{\vartheta}_k|^2_{\widehat{J}(k),\kappa_k}.
\end{aligned}
\tag{3.13}
$$

Using similar arguments, we have

$$
\begin{aligned}
|\vartheta_k|^2_{J(k)} = & \sum_{K\in\mathcal{M}_k} \sum_{\sigma\in\mathcal{F}_K} \frac{1}{d_{K,\sigma}} \int_\sigma [\vartheta_k]_{K,\sigma}(\boldsymbol{y})^2 d\boldsymbol{y} \\
\leq & \sum_{K\in\widehat{\mathcal{M}}_k} \sum_{\sigma\in\mathcal{F}_K} \frac{1}{d_{K,\sigma}} \int_\sigma [\hat{\vartheta}_k]_{K,\sigma,\kappa_k}(\boldsymbol{y})^2 d\boldsymbol{y} \\
& + 2 \sum_{T\in\mathcal{M}_k\setminus\widehat{\mathcal{M}}_k} \sum_{\sigma\in\mathcal{F}_T} \frac{1}{d_{T,\sigma}} \int_\sigma [\vartheta_{k-1}]_{T,\sigma}(\boldsymbol{y})^2 d\boldsymbol{y} \\
& + 2 \sum_{\substack{K\in\widehat{\mathcal{M}}_k \\ T\in\mathcal{M}_k\setminus\widehat{\mathcal{M}}_k}} \sum_{\sigma\in\mathcal{F}_K\cap\mathcal{F}_T} \frac{1}{d_{T,\sigma}} \int_\sigma \omega_{T,\sigma}^2 (\Pi_{\overline{K}}\hat{\vartheta}_k - \Pi_{\overline{K}}\vartheta_{k-1})(\boldsymbol{y})^2 d\boldsymbol{y} \\
\leq & |\hat{\vartheta}_k|^2_{\widehat{J}(k),\kappa_k} + 2|\vartheta_{k-1}|^2_{J(k)} + 2|\hat{\vartheta}_k - \vartheta_{k-1}|^2_{\partial\Omega_k^+}.
\end{aligned}
\tag{3.14}
$$

Combining Eqs. (3.13) and (3.14) we get

$$\|\nabla u - \nabla_{\mathcal{D}_k}\vartheta_k\|_{L^2(\Omega)^d} + |\vartheta_k|_{J(k)}$$
$$\leq \sqrt{2}(\|\nabla u - \nabla_{\mathcal{D}_{k-1}}\vartheta_{k-1}\|_{L^2(\Omega)^d} + \|\nabla u - \nabla_{\widehat{\mathcal{D}}_k,\kappa_k}\hat{\vartheta}_k\|_{L^2(\Omega_k)^d} + C_\psi|\hat{\vartheta}_k - \vartheta_{k-1}|_{\partial\Omega_k^+}$$
$$+ C_\psi|\hat{\vartheta}_k|_{\widehat{J}(k),\kappa_k}) + \sqrt{2}(|\hat{\vartheta}_k|_{\widehat{J}(k),\kappa_k} + |\vartheta_{k-1}|_{J(k)} + |\hat{\vartheta}_k - \vartheta_{k-1}|_{\partial\Omega_k^+})$$
$$\leq \sqrt{2}(\|\nabla u - \nabla_{\mathcal{D}_{k-1}}\vartheta_{k-1}\|_{L^2(\Omega)^d} + |\vartheta_{k-1}|_{J(k)} + \|\nabla u - \nabla_{\widehat{\mathcal{D}}_k,\kappa_k}\hat{\vartheta}_k\|_{L^2(\Omega_k)^d}$$
$$+ (1 + C_\psi)|\hat{\vartheta}_k|_{\widehat{J}(k),\kappa_k}) + \sqrt{2}(1 + C_\psi)|\hat{\vartheta}_k - \vartheta_{k-1}|_{\partial\Omega_k^+}.$$

Since we easily get

$$|\hat{\vartheta}_k - \vartheta_{k-1}|_{\partial\Omega_k^+} \leq \sqrt{2}C_{\omega,k}(|\hat{\vartheta}_k|_{\widehat{J}(k),\kappa_k} + |\vartheta_{k-1}|_{J(k)})$$

we obtain

$$\|\nabla u - \nabla_{\mathcal{D}_k}\vartheta_k\|_{L^2(\Omega)^d} + |\vartheta_k|_{J(k)}$$
$$\leq \sqrt{2}(1 + C_\psi)(1 + \sqrt{2}C_{\omega,k})(\|\nabla u - \nabla_{\mathcal{D}_k}\vartheta_{k-1}\|_{L^2(\Omega)^d} + |\vartheta_{k-1}|_{J(k)})$$
$$+ \sqrt{2}(1 + C_\psi)(1 + \sqrt{2}C_{\omega,k})(\|\nabla u - \nabla_{\widehat{\mathcal{D}}_k,\kappa_k}\hat{\vartheta}_k\|_{L^2(\Omega_k)^d} + |\hat{\vartheta}_k|_{\widehat{J}(k),\kappa_k}),$$

and the Lemma is proved. □

**Theorem 3.12.** *Let $u \in H_0^1(\Omega)$ be solution of Eq. (2.2) and $(\vartheta_k)_{k=1}^M$ be the sequence defined by scheme Eqs. (3.1a) to (3.1c). Then for $k = 1, ..., M$*

$$\lim_{h \to 0} \|\nabla u - \nabla_{\mathcal{D}_k}\vartheta_k\|_{L^2(\Omega)^d} + |\vartheta_k|_{J(k)} = 0. \tag{3.15}$$

*If moreover $u \in H_0^1(\Omega) \cap H^2(\Omega)$, the coefficients of A are Lipschitz continuous and $\boldsymbol{F} \in H^1(\Omega)^d$ there exists $C_6$ depending on $\alpha$, $\ell$, $d$, $\rho$, $C_r$, $|\Omega|$, $A$, $\boldsymbol{F}$ and $u$ such that*

$$\|\nabla u - \nabla_{\mathcal{D}_k}\vartheta_k\|_{L^2(\Omega)^d} + |\vartheta_k|_{J(k)} \leq C_6 h. \tag{3.16}$$

*Proof.* Follows from a recursive argument, Theorem 3.6 and Lemma 3.11. □

## 4 A priori error analysis for quasilinear problems

In this section we analyze our local SWDGGD scheme for a class of non linear problems satisfying Assumption 2.2. For the sake of simplicity we consider $f \in L^2(\Omega)$, but the algorithm and the results can easily be generalized to $f \in H^{-1}(\Omega)$. Under Assumption 2.2 there exists a unique weak solution $u \in H_0^1(\Omega)$ of

$$\int_\Omega A(u)\nabla u \cdot \nabla v \, d\boldsymbol{x} = \int_\Omega fv \, d\boldsymbol{x} \qquad \text{for all } v \in H_0^1(\Omega). \tag{4.1}$$

The local elliptic scheme for problem Eq. (4.1) is given as follows. Set $\vartheta_1 \in X_{\mathcal{D}_1}$ a solution of

$$\int_\Omega A(\Pi_{\mathcal{D}_1}\vartheta_1)\nabla_{\mathcal{D}_1}\vartheta_1 \cdot \nabla_{\mathcal{D}_1}\phi_1 \, d\boldsymbol{x} = \int_\Omega f\Pi_{\mathcal{D}_1}\phi_1 \, d\boldsymbol{x} \tag{4.2a}$$

for all $\phi_1 \in X_{\mathcal{D}_1}$. For $k \geq 2$ we set

$$\vartheta_k = \kappa_k + \hat{\vartheta}_k, \tag{4.2b}$$

where $\kappa_k \in Z_{\mathcal{D}_k}$ is given by

$$\kappa_k = \vartheta_{k-1}\chi_{\Omega\setminus\Omega_k} \tag{4.2c}$$

and $\hat{\vartheta}_k \in Y_{\mathcal{D}_k}$ is solution of

$$\int_{\Omega_k} A(\Pi_{\widehat{\mathcal{D}}_{k-1}}\hat{\vartheta}_{k-1})\nabla_{\widehat{\mathcal{D}}_k,\kappa_k}\hat{\vartheta}_k \cdot \nabla_{\widehat{\mathcal{D}}_k}\varphi_k \, d\boldsymbol{x} = \int_{\Omega_k} f\, \Pi_{\widehat{\mathcal{D}}_k}\varphi_k \, d\boldsymbol{x} \tag{4.2d}$$

for all $\varphi_k \in Y_{\mathcal{D}_k}$.

We define again a subset $\mathcal{H} \subset \mathbb{R}_+$ with zero as only accumulation point and for each $h \in \mathcal{H}$ two sequence of meshes $(\mathfrak{T}_{h,k})_{k=1}^M$, $(\widehat{\mathfrak{T}}_{h,k})_{k=1}^M$ satisfying Assumption 2.18 with $h = \max_{k=1,...,M} h_{\mathcal{M}_k}$. We consider the weighted gradient discretization methods $\mathcal{D}_{h,k}$, $\widehat{\mathcal{D}}_{h,k}$ deriving from $\mathfrak{T}_{h,k}$ and $\widehat{\mathfrak{T}}_{h,k}$, as defined in Section 2.3. The following theorem establishes the convergence of the non linear local SWDGGD scheme Eq. (4.2). The proof is inspired by a result in [9, chapter 2.1.4] for global non linear schemes.

**Theorem 4.1.** *For any $h \in \mathcal{H}$ there exists exactly one $\vartheta_{h,1} \in \mathcal{D}_{h,1}$ solution to Eq. (4.2a). Moreover, $\Pi_{\mathcal{D}_{h,1}}\vartheta_{h,1}$ converges strongly in $L^2(\Omega)$ to a solution $u$ of Eq. (4.1) and $\nabla_{\mathcal{D}_{h,1}}\vartheta_{h,1}$ converges strongly in $L^2(\Omega)^d$ to $\nabla u$ as $h \to 0$.*

We will prove that the same result holds for $\vartheta_{h,k}$ with $k \geq 2$. We start by proving convergence of the local solutions $\hat{\vartheta}_{h,k}$. For simplicity we drop the index $h$ in what follows.

**Theorem 4.2.** *Let Assumption 2.2 hold, $((\kappa_k, \hat{\vartheta}_k))_{k=1}^M$ be the sequence given by the local scheme Eqs. (4.2a) to (4.2d) and $u \in H_0^1(\Omega)$ be solution of Eq. (4.1). Then for $k = 1, ..., M$*

$$\lim_{h\to 0} \|\nabla u - \nabla_{\widehat{\mathcal{D}}_k,\kappa_k}\hat{\vartheta}_k\|_{L^2(\Omega_k)^d} = 0, \tag{4.3a}$$

$$\lim_{h\to 0} |\hat{\vartheta}_k|_{\widehat{J}(k),\kappa_k} = 0, \tag{4.3b}$$

*where the limit is taken for $h \in \mathcal{H}$.*

*Proof.* We will prove Eq. (4.3) by recursion. For $k = 1$ we easily get Eq. (4.3a), indeed $\kappa_1 = 0$, $\hat{\vartheta}_1 = \vartheta_1$ and by Theorem 4.1 we get

$$\lim_{h\to 0} \|\nabla u - \nabla_{\widehat{\mathcal{D}}_1,\kappa_1}\hat{\vartheta}_1\|_{L^2(\Omega_k)^d} = \lim_{h\to 0} \|\nabla u - \nabla_{\mathcal{D}_1}\vartheta_1\|_{L^2(\Omega)^d} = 0. \tag{4.4}$$

Let $\phi_1 \in X_{\mathcal{D}_1}$, we have

$$|\hat{\vartheta}_1|_{\widehat{J}(1),\kappa_1} = |\vartheta_1|_{J(1)} \leq |\vartheta_1 - \phi_1|_{J(1)} + |\phi_1|_{J(1)}.$$

From [11] we infer the existence of a constant $C_{\text{eq}}$ depending only on $\alpha$, $\ell$, $d$ such that

$$|\vartheta_1 - \phi_1|_{J(1)} \leq C_{\text{eq}}\|\nabla_{\mathcal{D}_1}\vartheta_1 - \nabla_{\mathcal{D}_1}\phi_1\|_{L^2(\Omega)^d},$$

hence
$$|\hat{\vartheta}_1|_{\widehat{J}(1),\kappa_1} \leq C_{\text{eq}} \|\nabla_{\mathcal{D}_1}\vartheta_1 - \nabla_{\mathcal{D}_1}\phi_1\|_{L^2(\Omega)^d} + |\phi_1|_{J(1)}.$$

Taking $\phi_1 = \operatorname{argmin}_{\phi \in X_{\mathcal{D}_1}}(\|\Pi_{\mathcal{D}_1}\phi - u\|_{L^2(\Omega)} + \|\nabla_{\mathcal{D}_1}\phi - \nabla u\|_{L^2(\Omega)^d} + |\phi|_{J(1)})$ we get Eq. (4.3b) for $k=1$ using the triangle inequality, Eq. (4.4) and Lemma 2.16.

Let $k \geq 2$ and suppose that Eq. (4.3) holds for $k-1$. By Lemma 3.1 there exists $\overline{\vartheta}_k \in Y_{\mathcal{D}_k}$ satisfying
$$\begin{aligned}\|\nabla_{\widehat{\mathcal{D}}_k,\kappa_k}\overline{\vartheta}_k - \nabla_{\widehat{\mathcal{D}}_{k-1},\kappa_{k-1}}\hat{\vartheta}_{k-1}\|_{L^2(\Omega_k)^d} &\leq C_{\text{i}}|\hat{\vartheta}_{k-1}|_{\widehat{J}(k-1),\kappa_{k-1}},\\ |\overline{\vartheta}_k|_{\widehat{J}(k),\kappa_k} &\leq C_{\text{i}}|\hat{\vartheta}_{k-1}|_{\widehat{J}(k-1),\kappa_{k-1}}.\end{aligned} \qquad (4.5)$$

Let $\tilde{\vartheta}_k \in Y_{\mathcal{D}_k}$ be solution of
$$\int_{\Omega_k} A_{k-1}(\nabla_{\widehat{\mathcal{D}}_k,\kappa_k}\overline{\vartheta}_k + \nabla_{\widehat{\mathcal{D}}_k}\tilde{\vartheta}_k)\cdot\nabla_{\widehat{\mathcal{D}}_k}\varphi_k\,d\boldsymbol{x} = \int_{\Omega_k} f\,\Pi_{\widehat{\mathcal{D}}_k}\varphi_k\,d\boldsymbol{x}$$
for all $\varphi_k \in Y_{\mathcal{D}_k}$, where $A_{k-1} = A(\Pi_{\widehat{\mathcal{D}}_{k-1}}\hat{\vartheta}_{k-1})$. Since $\nabla_{\widehat{\mathcal{D}}_k,\kappa_k}(\overline{\vartheta}_k + \tilde{\vartheta}_k) = \nabla_{\widehat{\mathcal{D}}_k,\kappa_k}\overline{\vartheta}_k + \nabla_{\widehat{\mathcal{D}}_k}\tilde{\vartheta}_k$ it follows that $\hat{\vartheta}_k = \overline{\vartheta}_k + \tilde{\vartheta}_k$. From Eq. (4.3) for $k-1$ and Eq. (4.5) it follows that $\nabla_{\widehat{\mathcal{D}}_k,\kappa_k}\overline{\vartheta}_k \to \nabla u$ strongly in $L^2(\Omega_k)^d$. Thus if $\nabla_{\widehat{\mathcal{D}}_k}\tilde{\vartheta}_k \to 0$ strongly in $L^2(\Omega_k)^d$ then $\nabla_{\widehat{\mathcal{D}}_k,\kappa_k}\hat{\vartheta}_k \to \nabla u$ strongly in $L^2(\Omega_k)^d$ and whence Eq. (4.3a) holds for $k$. From the coercivity of $A$
$$\begin{aligned}\underline{\lambda}\|\nabla_{\widehat{\mathcal{D}}_k}\tilde{\vartheta}_k\|^2_{L^2(\Omega_k)^d} &\leq \int_{\Omega_k} A_{k-1}\nabla_{\widehat{\mathcal{D}}_k}\tilde{\vartheta}_k\cdot\nabla_{\widehat{\mathcal{D}}_k}\tilde{\vartheta}_k\,d\boldsymbol{x}\\ &= \int_{\Omega_k} f\,\Pi_{\widehat{\mathcal{D}}_k}\tilde{\vartheta}_k\,d\boldsymbol{x} - \int_{\Omega_k} A_{k-1}\nabla_{\widehat{\mathcal{D}}_k,\kappa_k}\overline{\vartheta}_k\cdot\nabla_{\widehat{\mathcal{D}}_k}\tilde{\vartheta}_k\,d\boldsymbol{x}\\ &\leq \|f\|_{L^2(\Omega_k)}\|\Pi_{\widehat{\mathcal{D}}_k}\tilde{\vartheta}_k\|_{L^2(\Omega_k)} + \overline{\lambda}\|\nabla_{\widehat{\mathcal{D}}_k,\kappa_k}\overline{\vartheta}_k\|_{L^2(\Omega_k)^d}\|\nabla_{\widehat{\mathcal{D}}_k}\tilde{\vartheta}_k\|_{L^2(\Omega_k)^d}\\ &\leq (C_{\text{P}}\|f\|_{L^2(\Omega_k)} + \overline{\lambda}\|\nabla_{\widehat{\mathcal{D}}_k,\kappa_k}\overline{\vartheta}_k\|_{L^2(\Omega_k)^d})\|\nabla_{\widehat{\mathcal{D}}_k}\tilde{\vartheta}_k\|_{L^2(\Omega_k)^d}\end{aligned}$$

and hence $\|\nabla_{\widehat{\mathcal{D}}_k}\tilde{\vartheta}_k\|_{L^2(\Omega_k)^d}$ is bounded. It follows from the compactness of $\widehat{\mathcal{D}}_k$ and [9, Lemma 2.15] that there exists $w \in H^1_0(\Omega_k)$ and a subsequence $\mathcal{H}'$ of $\mathcal{H}$ such that $\Pi_{\widehat{\mathcal{D}}_k}\tilde{\vartheta}_k \to w$ strongly in $L^2(\Omega_k)$ and $\nabla_{\widehat{\mathcal{D}}_k}\tilde{\vartheta}_k \rightharpoonup \nabla w$ weakly in $L^2(\Omega_k)^d$ as $h \to 0$ with $h \in \mathcal{H}'$. We will show that $w = 0$, that the convergence holds for the whole sequence $\mathcal{H}$ and that $\nabla_{\widehat{\mathcal{D}}_k}\tilde{\vartheta}_k$ converges strongly. Let $v \in H^1_0(\Omega_k)$ and
$$\varphi_k = \operatorname{argmin}_{\varphi \in Y_{\mathcal{D}_k}}(\|\Pi_{\widehat{\mathcal{D}}_k}\varphi - v\|_{L^2(\Omega_k)} + \|\nabla_{\widehat{\mathcal{D}}_k}\varphi - \nabla v\|_{L^2(\Omega_k)^d} + |\varphi|_{\widehat{J}(k),0}).$$

Since $\widehat{\mathcal{D}}_k$ is a SWDGGD, from Lemma 2.16 we have that $\Pi_{\widehat{\mathcal{D}}_k}\varphi_k \to v$ strongly in $L^2(\Omega_k)$ and $\nabla_{\widehat{\mathcal{D}}_k}\varphi_k \to \nabla v$ strongly in $L^2(\Omega_k)^d$. From Eq. (4.3a) $\nabla_{\widehat{\mathcal{D}}_{k-1},\kappa_{k-1}}\hat{\vartheta}_{k-1} \to \nabla u$ strongly in $L^2(\Omega_k)^d$, furthermore by coercivity and consistency we can show that $\Pi_{\widehat{\mathcal{D}}_{k-1}}\hat{\vartheta}_{k-1} \to u$ strongly in $L^2(\Omega_k)$ as well. The same holds for $\overline{\vartheta}_k$. Hence by the non-linear strong convergence Lemma [9, section D.4] we obtain
$$\begin{aligned}A(\Pi_{\widehat{\mathcal{D}}_{k-1}}\hat{\vartheta}_{k-1})\nabla_{\widehat{\mathcal{D}}_k}\varphi_k &\to A(u)\nabla v \text{ strongly in } L^2(\Omega_k)^d,\\ A(\Pi_{\widehat{\mathcal{D}}_{k-1}}\hat{\vartheta}_{k-1})\nabla_{\widehat{\mathcal{D}}_k,\kappa_k}\overline{\vartheta}_k &\to A(u)\nabla u \text{ strongly in } L^2(\Omega_k)^d.\end{aligned}$$

It follows from weak-strong convergence Lemma [9, section D.4] and symmetry of $A$ that

$$\int_{\Omega_k} A(u)\nabla w \cdot \nabla v \, d\boldsymbol{x} = \int_{\Omega_k} \nabla w \cdot A(u)\nabla v \, d\boldsymbol{x}$$
$$= \lim_{h\to 0} \int_{\Omega_k} \nabla_{\widehat{\mathcal{D}}_k} \tilde{\vartheta}_k \cdot A_{k-1} \nabla_{\widehat{\mathcal{D}}_k} \varphi_k \, d\boldsymbol{x}, \qquad (4.6)$$

where the limit is for $h \in \mathcal{H}'$. On the other hand we have

$$\int_{\Omega_k} A_{k-1} \nabla_{\widehat{\mathcal{D}}_k} \tilde{\vartheta}_k \cdot \nabla_{\widehat{\mathcal{D}}_k} \varphi_k \, d\boldsymbol{x} = \int_{\Omega_k} f \, \Pi_{\widehat{\mathcal{D}}_k} \varphi_k \, d\boldsymbol{x}$$
$$- \int_{\Omega_k} A_{k-1} \nabla_{\widehat{\mathcal{D}}_k, \kappa_k} \overline{\vartheta}_k \cdot \nabla_{\widehat{\mathcal{D}}_k} \varphi_k \, d\boldsymbol{x}$$

and taking the limit we get

$$\lim_{h\to 0} \int_{\Omega_k} A_{k-1} \nabla_{\widehat{\mathcal{D}}_k} \tilde{\vartheta}_k \cdot \nabla_{\widehat{\mathcal{D}}_k} \varphi_k \, d\boldsymbol{x} = \int_{\Omega_k} fv \, d\boldsymbol{x} - \int_{\Omega_k} A(u)\nabla u \cdot \nabla v \, d\boldsymbol{x} = 0. \qquad (4.7)$$

Putting toghether Eqs. (4.6) and (4.7) and using the symmetry of $A_{k-1}$ we obtain

$$\int_{\Omega_k} A(u)\nabla w \cdot \nabla v \, d\boldsymbol{x} = 0$$

for all $v \in H_0^1(\Omega_k)$ and so $w = 0$. We can repeat the same reasoning for each subsequence of $\nabla_{\widehat{\mathcal{D}}_k} \tilde{\vartheta}_k$ and obtain the same limit $w = 0$, hence $\Pi_{\widehat{\mathcal{D}}_k} \tilde{\vartheta}_k \to 0$ strongly in $L^2(\Omega_k)$ and $\nabla_{\widehat{\mathcal{D}}_k} \tilde{\vartheta}_k \rightharpoonup 0$ weakly in $L^2(\Omega_k)^d$ for the whole sequence $\mathcal{H}$. Furthermore

$$\int_{\Omega_k} A_{k-1} \nabla_{\widehat{\mathcal{D}}_k} \tilde{\vartheta}_k \cdot \nabla_{\widehat{\mathcal{D}}_k} \tilde{\vartheta}_k \, d\boldsymbol{x} = \int_{\Omega_k} f \, \Pi_{\widehat{\mathcal{D}}_k} \tilde{\vartheta}_k \, d\boldsymbol{x}$$
$$- \int_{\Omega_k} A_{k-1} \nabla_{\widehat{\mathcal{D}}_k, \kappa_k} \overline{\vartheta}_k \cdot \nabla_{\widehat{\mathcal{D}}_k} \tilde{\vartheta}_k \, d\boldsymbol{x}$$

and so

$$\lim_{h\to 0} \int_{\Omega_k} A_{k-1} \nabla_{\widehat{\mathcal{D}}_k} \tilde{\vartheta}_k \cdot \nabla_{\widehat{\mathcal{D}}_k} \tilde{\vartheta}_k \, d\boldsymbol{x} = 0,$$

which shows that $\lim_{h\to 0} \|\nabla_{\widehat{\mathcal{D}}_k} \tilde{\vartheta}_k\|_{L^2(\Omega_k)^d} = 0$ and hence the strong convergence of $\nabla_{\widehat{\mathcal{D}}_k} \tilde{\vartheta}_k$. It rests to show Eq. (4.3b), we have

$$|\hat{\vartheta}_k|_{\widehat{J}(k),\kappa_k} \leq |\overline{\vartheta}_k|_{\widehat{J}(k),\kappa_k} + |\tilde{\vartheta}_k|_{\widehat{J}(k),0}$$
$$\leq C_{\mathrm{i}} |\vartheta_{k-1}|_{\widehat{J}(k-1),\kappa_{k-1}} + C_{\mathrm{eq}} \|\nabla_{\widehat{\mathcal{D}}_k} \tilde{\vartheta}_k\|_{L^2(\Omega_k)^d}$$

and the result follows. $\square$

The next Theorem can be proved with similar arguments as used in Section 3.

**Theorem 4.3.** *Let Assumption 2.2 hold. Consider $(\vartheta_k)_{k=1}^M$, the sequence given by the local scheme Eqs. (4.2a) to (4.2d) and $u \in H_0^1(\Omega)$ the solution of Eq. (4.1). Then for $k = 1, ..., M$, we have*

$$\lim_{h\to 0} \|\nabla u - \nabla_{\mathcal{D}_k} \vartheta_k\|_{L^2(\Omega)^d} + |\vartheta_k|_{J(k)} = 0,$$

*where the limit is taken for $h \in \mathcal{H}$.*

# 5 Numerical Experiments

In the following numerical experiments, we will use examples where the subdomains $\{\Omega_k\}_{k=1}^M$ and meshes $\{\mathfrak{T}_k\}_{k=1}^M$ are defined a priori. This might be realistic in applications where the location of the singularities or high contrast of the solution are known a priori. When such a priori knowledge is not available, we should use instead a posteriori error estimators for detecting the local subdomains. This is developed in a companion paper [1].

In what follows $\{\Omega_k\}_{k=1}^M$ will be a sequence of embedded domains but we recall that this is not a requirement. In the examples we consider $f \in L^2(\Omega)$ and denote by $\zeta_k \in X_{\mathcal{D}_k}$ the solution of

$$\int_\Omega A(\Pi_{\mathcal{D}_k}\zeta_k)\nabla_{\mathcal{D}_k}\zeta_k \cdot \nabla_{\mathcal{D}_k}\phi_k\, d\boldsymbol{x} = \int_\Omega f\, \Pi_{\mathcal{D}_k}\phi_k\, d\boldsymbol{x} \qquad \text{for all } \phi_k \in X_{\mathcal{D}_k}, \tag{5.1}$$

we refer to $\zeta_k$ as the classical solution, that is, the one obtained by the usual scheme which solves the equations in the whole domain after each mesh refinement. We can write $\zeta_k = \hat{\zeta}_k + \eta_k$ with $\hat{\zeta}_k \in Y_{\mathcal{D}_k}$ and $\eta_k \in Z_{\mathcal{D}_k}$. We will often compare $\vartheta_k$ and $\hat{\vartheta}_k$ the solutions of Eq. (3.1) or Eq. (4.2) against $\zeta_k$ and $\hat{\zeta}_k$ respectively.

**Computational cost.** As in our setting, the meshes are defined a priori only the the most accurate solution $\zeta_M$ need to be computed. For the iterative schemes Eqs. (3.1) and (4.2) instead it is imperative to compute $\vartheta_k$ for $k = 1, ..., M$. If for example a conjugate gradient method is used to solve the linear systems, then the computational cost of the local scheme can be considerably smaller than the classical scheme due to the smaller problems solved on the fine meshes. For nonlinear problems, the local scheme might be faster for any linear solver, as the non linear system is solved only on a coarse mesh (see Section 4). This is illustrated in our numerical experiments.

It is useful to define the quantities

$$\text{Local Err}(\hat{\vartheta}_k) := \|\nabla u - \nabla_{\widehat{\mathcal{D}}_k,\kappa_k}\hat{\vartheta}_k\|_{L^2(\Omega_k)^d} + |\hat{\vartheta}_k|_{\widehat{\mathcal{J}}(k),\kappa_k},$$
$$\text{Global Err}(\vartheta_k) := \|\nabla u - \nabla_{\mathcal{D}_k}\vartheta_k\|_{L^2(\Omega)^d} + |\vartheta_k|_{\mathcal{J}(k)}.$$

Similarly we define Local Err$(\hat{\zeta}_k)$ and Global Err$(\zeta_k)$ for the local and global error of the classical solutions. The local and cassical schemes have been implemented with the help of the C++ library `libMesh` [16].

## 5.1 Convergence Rates

In this example we want to verify results Eqs. (3.16), (3.7b) and (3.7c), hence we consider an example with smooth solution. Let $\Omega = [-1, 1] \times [-1, 1]$, $A = I_2$ the identity matrix and $f \in L^2(\Omega)$ such that the exact solution is

$$u(\boldsymbol{x}) = e^{-120\|\boldsymbol{x}\|_2^2}. \tag{5.2}$$

Let $M = 4$, the local domains are such that $\boldsymbol{x} \in \Omega_k$ if $\|\boldsymbol{x}\|_\infty < (5-k)/4$ for $k = 1, ..., 4$.

In the first experiment we want to verify the estimates Eqs. (3.16) and (3.7b), i.e. the convergence of the local and global errors with respect to the global mesh size. For a fixed $h$ we consider uniform simplicial meshes $\widehat{\mathcal{M}}_k$ on $\Omega_k$ with mesh size $h_{\widehat{\mathcal{M}}_k} = h/2^{k-1}$ and apply the local algorithm Eq. (3.1), we let $h \to 0$ and verify the convergence rates. From Figs. 3(a) and 3(b) we see that Eqs. (3.16) and (3.7b) are verified for the local solution $\vartheta_4$. We also see that the classical scheme gives results with the same accuracy as the local scheme, both for the local and the global error. This example also indicates that if the high gradient regions are localized then there is no need of solving the problem in the whole domain after refinements.

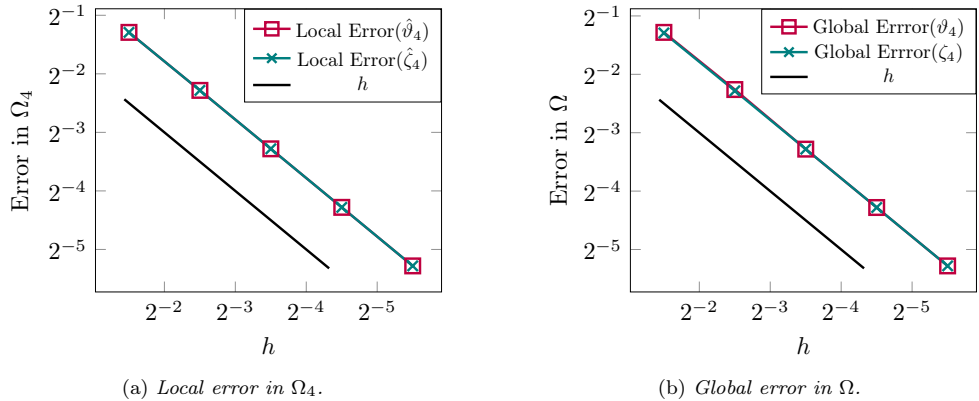

(a) *Local error in $\Omega_4$.*

(b) *Global error in $\Omega$.*

Figure 3. *Experiment 5.1: convergence of the local $\vartheta_4$ and classical $\zeta_4$ solutions letting $h \to 0$.*

In the next experiment we want to see the influence of the second term (boundary layer term) in the righthand side of Eq. (3.7c) on Local Error($\hat{\vartheta}_M$). Let $r \in\, ]0, 1[$, we set $M = 2$, $\Omega_1 = \Omega$ and $\Omega_2 = [-r, r] \times [-r, r]$. We fix $h_{\mathcal{M}_1} = \sqrt{2}/8$ the mesh size of $\mathcal{M}_1$ and let $h_{\widehat{\mathcal{M}}_2} \to 0$. We plot the results for different values of $r$ (an illustration of this numerical experiment is given in Fig. 5). In

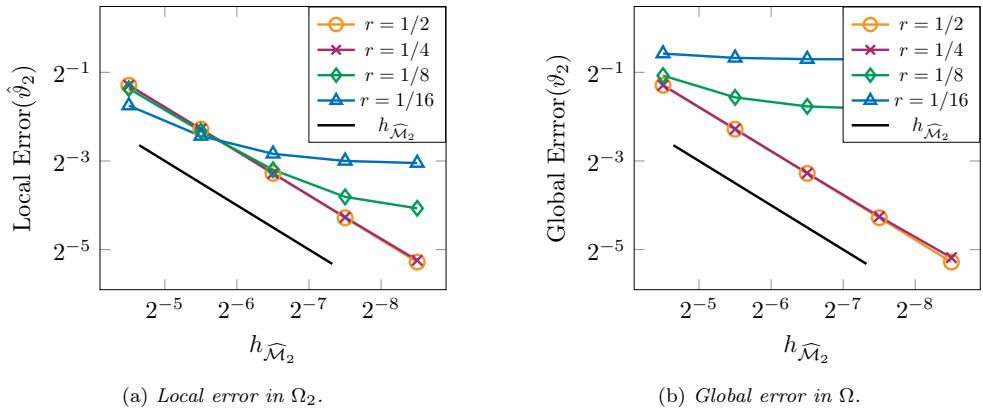

(a) *Local error in $\Omega_2$.*

(b) *Global error in $\Omega$.*

Figure 4. *Experiment 5.1: effect of the size of $\Omega_2$ on the local solution $\vartheta_2$ when $h_{\widehat{\mathcal{M}}_2} \to 0$.*

Fig. 4(a) we see that when $r$ is large enough the local error scales with the local mesh size. If, instead, $r$ is too small to cover the high gradient regions then the local error saturates very quickly. With $r = 1/8$ we get nice convergence up to $h_{\mathcal{M}_1}/h_{\widehat{\mathcal{M}}_2} = 16$ and with $r = 1/4, 1/2$ we do not see any saturation effects. In Fig. 5 we see that $r = 1/16$ is too small to cover the local variations and indeed the local error does not converges. In Fig. 4(b) we plot the total error on $\Omega$. We remark that the error saturates for $r = 1/16, 1/8$. It is interesting to compare the results for $r = 1/8$ in Figs. 4(a)

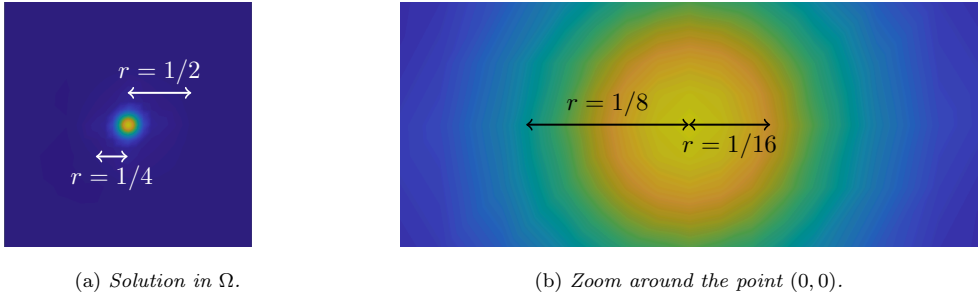

(a) *Solution in $\Omega$.*    (b) *Zoom around the point $(0,0)$.*

Figure 5. *Experiment 5.1: solution $u$ from Eq. (5.2) with the size of domains $\Omega_2$ depending on $r$.*

and 4(b), in the first one there is a nice convergence while in the second an immediate saturation. This indicates that even if the error outside of $\Omega_2$ is important, it does not propagate quickly into $\Omega_2$.

In Fig. 6 we plot the results of the same experiment shown in Fig. 4 but for $\zeta_4$ instead of $\vartheta_4$. We see that, again, the classical scheme gives similar results.

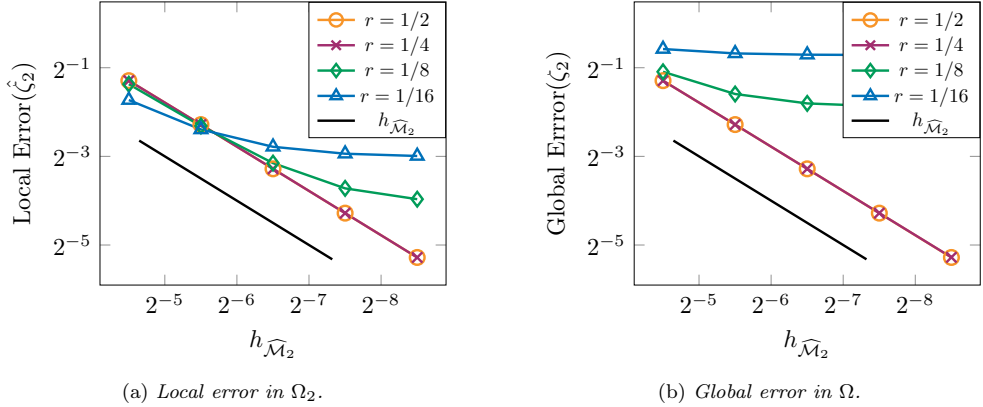

(a) *Local error in $\Omega_2$.*    (b) *Global error in $\Omega$.*

Figure 6. *Experiment 5.1: effect of the size of $\Omega_2$ on the classical solution $\zeta_2$ when $h_{\widehat{\mathcal{M}}_2} \to 0$.*

## 5.2 Influence of artificial boundary conditions

The goal of this experiment is to verify the result of Theorem 3.9, we want to illustrate numerically that the error due to artificial boundary conditions is of higher order as proved in estimate Eq. (3.12). We consider the same problem as in Section 5.1 with $M = 2$, $\Omega_1 = \Omega$ and $\Omega_2 = [-r, r] \times [-r, r]$ with $r = 1/16$. We saw previously that with this choice of $r$ the error originating from the artificial boundary conditions dominates the local error in $\Omega_2$. We solve Eq. (3.1) with different mesh sizes $h = h_{\mathcal{M}_1}$ using $h_{\widehat{\mathcal{M}}_2} = h/2^5$ as local mesh size, with this choice of $h_{\widehat{\mathcal{M}}_2}$ the dominating term in Eq. (3.12) is the last one, i.e. the one in $h^{3/2} \log(h^{-1})$. We measure the local errors in $\Omega_2$ and plot the results in Fig. 7. We see that indeed the local error satisfies Eq. (3.12) and converges even slightly faster than predicted.

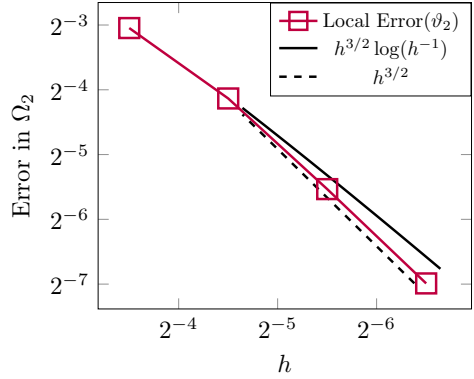

Figure 7. *Experiment 5.2: convergence order of artificial boundary conditions error term.*

## 5.3 Non regular problem: discontinuous data

We next want to probe numerically the convergence of our local scheme for a solution only belonging to $H^{1+\varepsilon}(\Omega)$ (for small $\varepsilon > 0$). The convergence is predicted by estimates Eqs. (3.15) and (3.7a). We consider a problem that has been studied in [15] and [20] (in the context of a posteriori error estimators).

Let $\Omega = [-1,1] \times [-1,1]$ and consider Problem Eq. (1.1) with $f = 0$. We divide the computational domain in four equal parts. Let the tensor be defined as $A(\boldsymbol{x}) = a_1 I_2$ in the 1st and 3rd quadrants and $A(\boldsymbol{x}) = a_2 I_2$ in the 2nd and 4th quadrants. The exact solution is given by $u(r,\theta) = r^\gamma \mu(\theta)$, where

$$\mu(\theta) = \begin{cases} \cos((\pi/2 - \sigma)\gamma) \cos((\theta - \pi/2 + \rho)\gamma) & \text{if } 0 \leq \theta \leq \pi/2, \\ \cos(\rho\gamma) \cos((\theta - \pi + \sigma)\gamma) & \text{if } \pi/2 < \theta \leq \pi, \\ \cos(\sigma\gamma) \cos((\theta - \pi - \rho)\gamma) & \text{if } \pi < \theta \leq 3\pi/2, \\ \cos((\pi/2 - \rho)\gamma) \cos((\theta - 3\pi/2 - \sigma)\gamma) & \text{if } 3\pi/2 < \theta < 2\pi. \end{cases}$$

The parameters $\gamma$, $\rho$, $\sigma$ and $R := a_1/a_2$ satisfy the following non linear equations

$$R = -\tan((\pi/2 - \sigma)\gamma) \cot(\rho\gamma),$$
$$1/R = -\tan(\rho\gamma) \cot(\sigma\rho),$$
$$R = -\tan(\rho\gamma) \cot((\pi/2 - \rho)\gamma),$$
$$\max\{0, \pi\gamma - \pi\} < 2\gamma\rho < \min\{\pi\gamma, \pi\},$$
$$\max\{0, \pi - \pi\gamma\} < -2\gamma\sigma < \min\{\pi, 2\pi - \pi\gamma\}.$$

It is known that $u \in H^{1+\gamma-\epsilon}(\Omega)$ for any $\epsilon > 0$. In this example we choose $\gamma = 0.1$, $\sigma = -19\pi/4$, $\rho = \pi/4$ and $R \approx 161$.

In order to verify the estimates Eqs. (3.15) and (3.7a), we perform the same experiments as in Section 5.1, shown in Fig. 3. We take $M = 4$ and the same domain and mesh sequences. We let $h \to 0$ and show the results in Fig. 8. We find a convergence rate of 0.09, which is consistent with the results of [17] and the fact that $u$ is almost in $H^{1.1}(\Omega)$. As was observed in Section 5.1, we see that the two solutions $\vartheta_4$ and $\zeta_4$ have the same errors, both in the local and global domains.

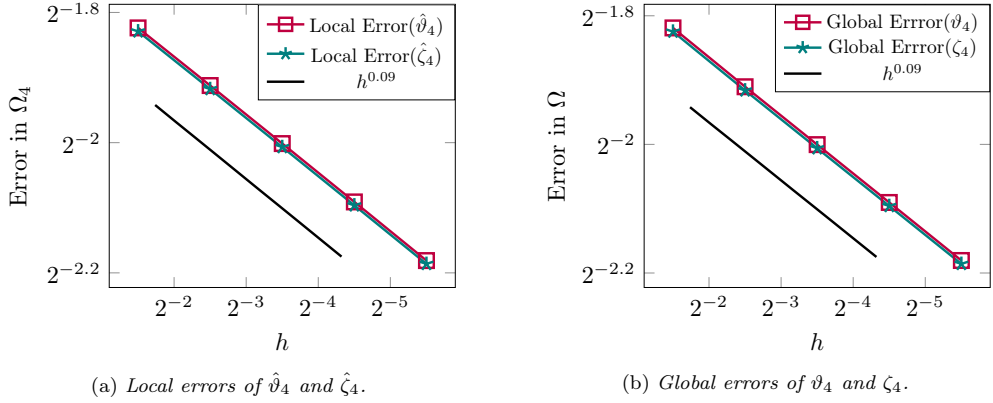

(a) *Local errors of $\hat{\vartheta}_4$ and $\hat{\zeta}_4$.*

(b) *Global errors of $\vartheta_4$ and $\zeta_4$.*

Figure 8. *Experiment 5.3: convergence of $\zeta_4$ and $\vartheta_4$ letting $h \to 0$.*

The influence of the term $|\kappa_k - \xi_k|_{\partial\Omega_k^-}$ in Eq. (3.3) on Local Error$(\hat{\vartheta}_M)$ is established next, repeating the experiment of Section 5.1, taking $\Omega_2$ depending on $r \in ]0,1[$ and letting $h_{\widehat{\mathcal{M}}_2} \to 0$. The results for $\vartheta_2$ and $\zeta_2$ are plotted in Figs. 9 and 10 respectively. In contrast to the previous experiment, we do not have any saturation since the error inside the local domain largely dominates.

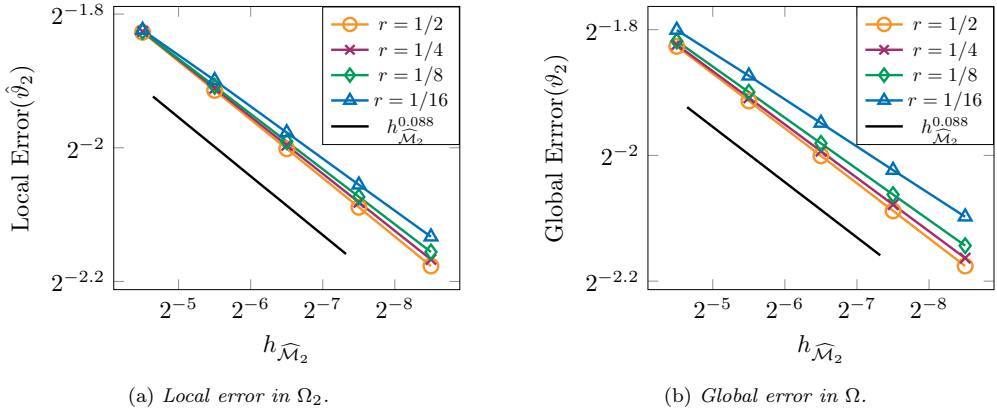

(a) *Local error in $\Omega_2$.*

(b) *Global error in $\Omega$.*

Figure 9. *Experiment 5.3: effect of the size of $\Omega_2$ on the local solution $\vartheta_2$ when $h_{\widehat{\mathcal{M}}_2} \to 0$.*

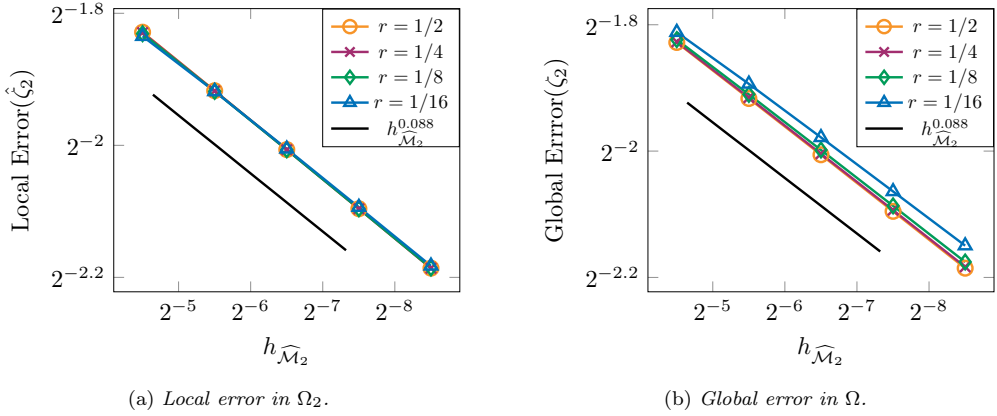

(a) *Local error in $\Omega_2$.*   (b) *Global error in $\Omega$.*

Figure 10. *Experiment 5.3: effect of the size of $\Omega_2$ on the classical solution $\zeta_2$ when $h_{\widehat{\mathcal{M}}_2} \to 0$.*

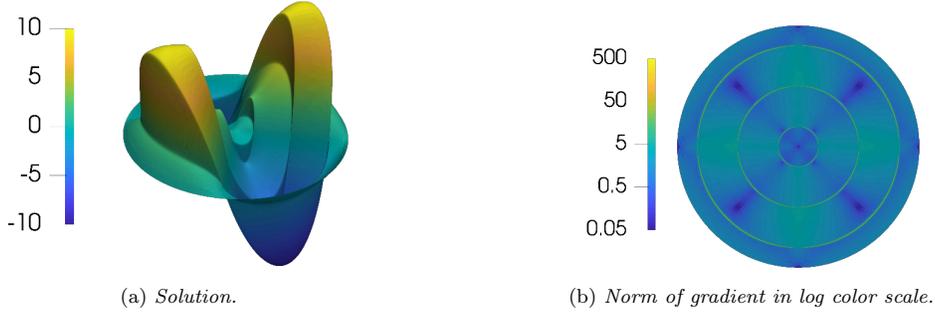

(a) *Solution.*   (b) *Norm of gradient in log color scale.*

Figure 11. *Experiment 5.4: solution and norm of the gradient.*

## 5.4 Computational cost of local versus non-local scheme for a linear equation

In this experiment we want to compare the numerical efficiency of the classical and local schemes on a linear equation, by computing a sequence of solutions with each scheme and plotting the accuracy against the cost.

We consider equation Eq. (1.1) with $\Omega = \{x \in \mathbb{R}^2 : \|x\|_2 < 3\pi\}$, a diffusion tensor $A(x) = \varepsilon + 1 - \sin(\|x\|_2)^{100}$ with $\varepsilon = 10^{-3}$ and the force $f$ is 1 if $x$ is the first or third quadrants and $-1$ else. An illustration of the solution is given in Fig. 11. We choose five local domains defined as $\Omega_1 = \Omega$ and

$$\Omega_k = \bigcup_{j=1}^{3} \{x \in \mathbb{R}^2 : |\|x\|_2 - (2j-1)\pi/2| < 2^{2-k}\} \quad \text{for } k = 2, ..., 5.$$

The meshes $\widehat{\mathcal{M}}_k$ are built so that $h_{\widehat{\mathcal{M}}_1} \approx 0.3$ and for $k = 2, ..., 5$ we have $h_{\widehat{\mathcal{M}}_k} = h_{\widehat{\mathcal{M}}_{k-1}}/2$. We run the local scheme and at each level we compute the full error and cost of $\vartheta_k$. As a measure of the cost for $\vartheta_k$ we take the sum of the time spent solving the linear systems up to level $k$ using the conjugate gradient (CG) method with incomplete Cholesky (IC) factorization as preconditioner. In [13] it is shown that this approach is the most robust and efficient for such problems. Then we run the classical scheme Eq. (5.1) on each mesh $\mathcal{M}_k$ and obtain a sequence of solutions $\zeta_k$. For each $k = 1, ..., 5$ we compute the full error and cost of $\zeta_k$. The cost is given by the time spent for solving the linear system at level $k$, where we use again CG with IC as preconditioner. Observe that here the cost is not cumulative as in the local method, since the classical scheme does not need $\zeta_{k-1}$ in order to compute $\zeta_k$. In Fig. 12(a) we plot the global error against the cost for both schemes, we see a significant speedup for the local scheme. In Fig. 12(b) we plot the speed up in function of the error, the graph is obtained dividing the two curves seen in Fig. 12(a).

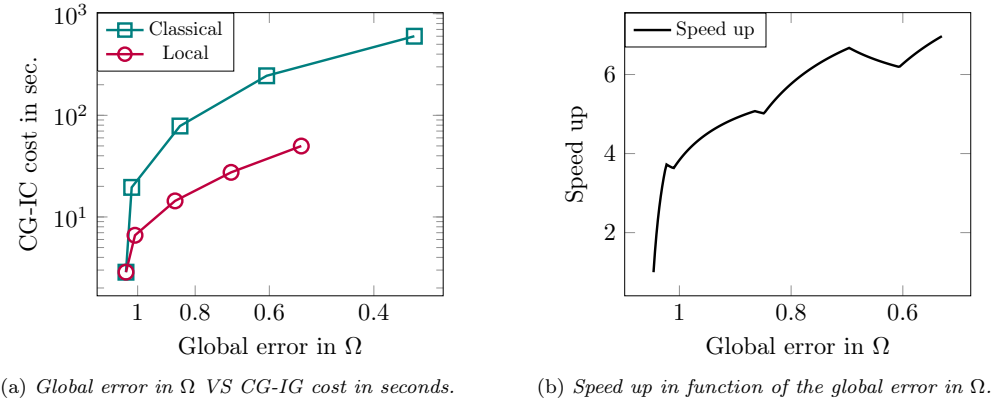

(a) *Global error in $\Omega$ VS CG-IG cost in seconds.*   (b) *Speed up in function of the global error in $\Omega$.*

Figure 12. *Experiment 5.4: performance comparison in a linear case.*

For linear problems such as in this experiment, the reason for the speed up is not only the reduced number of degrees of freedom but mostly the condition number of the linear system. The classical scheme solves linear systems arising from FE discretizazion on the whole domain, hence the matrix has high variations in its components due to possibly high contrasts in the tensor and the variation in the measure of the different elements. Instead, the local scheme uses matrices built from local discretizations, hence the tensor has milder variations and the elements of the local mesh have uniform size. This leads to matrices with smaller condition number. We see in Fig. 13(a) that the number of degrees of freedom of the two schemes is almost the same, while in Fig. 13(b) it is shown that the condition number of the stiffness matrix is much lower for the local scheme.

## 5.5 Quasilinear equation

In our last numerical experiment we want to compare the efficiency of the local and classical methods when solving a quasilinear equation. We consider the stationary Richards equation in pressure head form, given by

$$-\nabla \cdot (A(\boldsymbol{x}, h)\nabla(h - x_2)) = 0. \quad (5.3)$$

It describes the movement of a fluid in an unsaturated media and can be put in the form of Eq. (4.1) with the change of variables $u = h - x_2$. We consider $\Omega = [-50, 50] \times [-50, 50]$ and add the Dirichlet condition $g(\boldsymbol{x}) = 10(50 - x_2) + 3(50 + x_2)$. The diffusion tensor is given by $A(\boldsymbol{x}, h) = A_s(\boldsymbol{x})A_r(h)$,

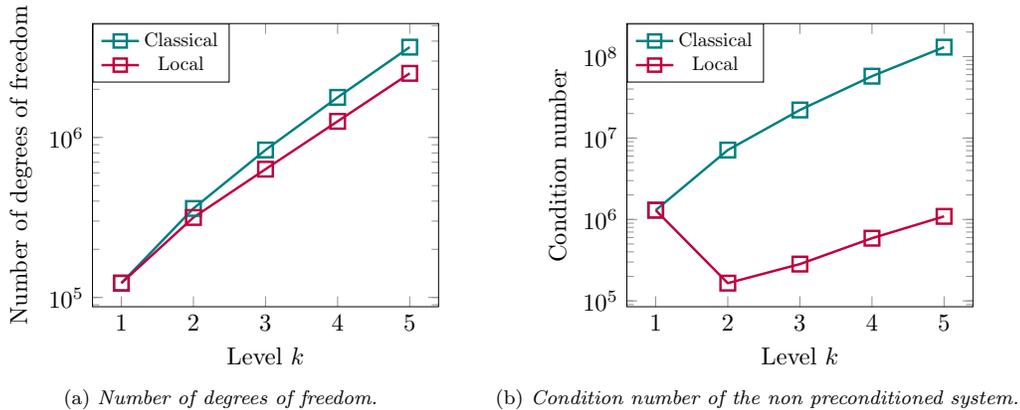

(a) Number of degrees of freedom.

(b) Condition number of the non preconditioned system.

Figure 13. *Experiment 5.4: properties of the linear systems.*

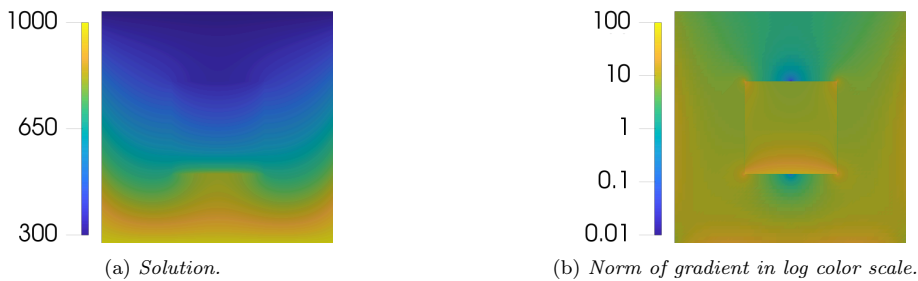

(a) Solution.

(b) Norm of gradient in log color scale.

Figure 14. *Experiment 5.5: solution and norm of the gradient for the Richards equation.*

where $A_s(\boldsymbol{x})$ is the conductivity in saturated conditions and $A_r(h)$ is the relative conductivity. These latter quantities are defined by

$$A_s(\boldsymbol{x}) = \begin{cases} 10^{-3} & \text{if } \|\boldsymbol{x}\|_\infty \leq 20, \\ 1 & \text{else,} \end{cases} \qquad A_r(h) = \frac{(1-(ah)^{n-1}(1+(ah)^n)^{-m})^2}{(1+(ah)^n)^{m/2}}.$$

The model $A_r(h)$ has been taken from [22], where $m = 1 - 1/n$ is chosen. The parameters $a, n$ are soil dependent: we choose $a = 1/500$ and $n = 2.68$, which is in the range of real case parameters. Remark that the tensor is discontinuous in $\boldsymbol{x}$ and hence does not satisfy Assumption 2.2. We plot in Fig. 14 the reference solution and the norm of its gradient, we see that the gradient is highly discontinuous.

Let $M = 4$, $\Omega_1 = \Omega$ and $\Omega_k$ for $k = 2, 3, 4$ defined by $\boldsymbol{x} \in \Omega_k$ if $\|\boldsymbol{x}\|_\infty \leq 20(1 + 2^{-k})$. First, we fix $h_{\widehat{\mathcal{M}}_k} = 100\sqrt{2}/2^{4+k}$ and compute the local solutions $\vartheta_k$ given by the local method Eq. (4.2). At the first level $k = 1$ we need to solve a nonlinear problem on a coarse grid using Newton iterations, where the initial guess is an extrapolation of the Dirichlet condition $g(\boldsymbol{x})$ on the whole domain. In the next levels $k > 1$ the local scheme solves a linear system using the Picard iteration step defined in Eq. (4.2d). At each level we compute the full error and cost of $\vartheta_k$. As a measure of the cost for $\vartheta_k$ we take the sum of the time spent solving the linear and non linear systems up to level $k$.

Since at $k = 1$ we perform a linearization of the system, it is no more symmetric because of the additional term, hence it has to be solved with the GMRES iterative scheme with incomplete LU (ILU) factorization as preconditioner, instead of CG with IC. In the following iterations with $k \geq 2$ we solve a linear system and hence the CG scheme with IC is used.

Then we compute similar solutions with the classical method and compare the costs. For the classical solution we need, for each $k = 1, 2, 3, 4$, to solve Eq. (5.3) with the Newton method. As initial guess we take again $g(\boldsymbol{x})$ and the Newton iterations are stopped when the error of the classical solution $\zeta_k$ is similar to the one of $\vartheta_k$. In about 3 or 4 Newton iterations we obtained errors differing by only about 1%. To measure the cost for $\zeta_k$ we consider the time spent in solving the non linear system at level $k$. The cost here is not cumulative as in the local method but on the other hand the linear systems to solve are not symmetric and the GMRES scheme with ILU preconditioner is used. In Fig. 15(a) we plot the error against the cost for this experiment. We see that the local scheme performs much better than the classical scheme in terms of computational cost versus accuracy.

Finally we compare the accuracy and cost of solving the local systems Eq. (4.2d) replacing $\Pi_{\widehat{\mathcal{D}}_{k-1}} \hat{\vartheta}_{k-1}$ with $\Pi_{\widehat{\mathcal{D}}_k} \hat{\vartheta}_k$, i.e., defining nonlinear local problems. These local systems have now to be solved by Newton method and GMRES with ILU. We denote by $\theta_k^1$ the solution where we use one Newton iteration and by $\theta_k^2$ the solution with two Newton iterations. In Fig. 15(b) we plot the error against the cost for $\vartheta_k$, $\theta_k^1$ and $\theta_k^2$. We see that one Picard iteration gives very similar error to the one or two Newton iterations but at a smaller cost, thanks to the CG scheme.

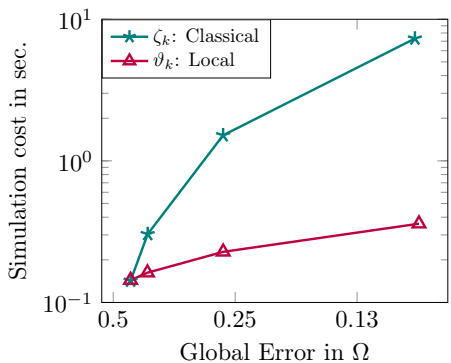

(a) *Efficiency of local and classical schemes.*

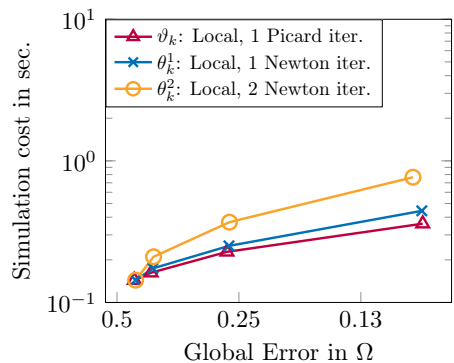

(b) *Efficiency of local scheme and two variants*

Figure 15. *Experiment 5.5: performance of classical scheme, local scheme and local scheme with Newton iterations instead of one Picard iteration.*

# 6 Conclusion

In this paper we introduced a local scheme for linear and quasilinear elliptic equations. The method does not rely on an iterative procedure and only needs one global solve on a coarse mesh. All subsequent computations are local. The a priori error analysis has been performed under weak regularity assumptions thanks to the gradient discretization framework. Numerical experiments have shown the efficiency of the scheme when applied to equations with localized high gradient regions. In a forthcoming paper [1] the a posteriori error analysis of the same scheme will be presented. Thanks to the a posteriori error estimators the local domains can be defined even when the singularities are not known a priori. We note that the extension of the local scheme to parabolic problems is also of interest. In particular, we believe that the techniques developed for a priori and a posteriori error analysis for the local scheme for elliptic PDEs also allow to analyze local time stepping schemes for parabolic PDEs.

## Acknowledgements

The authors are partially supported by the Swiss national foundation grant 200020_172710.

## A  Equivalence to SWIPG scheme

In this appendix we show that the SWDGGD scheme described in Section 2.2 is equivalent to the SWIPG method of [10, 7]. In particular we show that Eq. (2.3) with $\mathcal{D} = (X_\mathcal{D}, \Pi_\mathcal{D}, \nabla_\mathcal{D})$ as defined in Section 2.2 is equivalent to [7, equation 4.63], in order to do that we follow [11], where the equivalence of a GD to the SIP method is shown. We suppose that $A(\boldsymbol{x}, u) = A(\boldsymbol{x})$ and $A_K := A|_K$ the restriction of $A$ to an element $K \in \mathcal{M}$ is constant, that $f \in L^2(\Omega)$ and hence $f = f_0$.

Starting from Eq. (2.3) and developing the gradients we get

$$\int_\Omega A \nabla_\mathcal{D} \vartheta \cdot \nabla_\mathcal{D} \phi \, d\boldsymbol{x}$$

$$= \sum_{K \in \mathcal{M}} \int_K A_K \nabla_{\overline{K}} \vartheta \cdot \nabla_{\overline{K}} \phi \, d\boldsymbol{x}$$

$$+ \sum_{K \in \mathcal{M}} \sum_{\sigma \in \mathcal{F}_K} \int_{D_{K,\sigma}} \frac{\psi(s)}{d_{K,\sigma}} A_K ([\vartheta]_{K,\sigma}(\boldsymbol{y}) \nabla_{\overline{K}} \phi(\boldsymbol{x}) + [\phi]_{K,\sigma}(\boldsymbol{y}) \nabla_{\overline{K}} \vartheta(\boldsymbol{x})) \cdot \boldsymbol{n}_{K,\sigma} \, d\boldsymbol{x}$$

$$+ \sum_{K \in \mathcal{M}} \sum_{\sigma \in \mathcal{F}_K} \int_{D_{K,\sigma}} \frac{A_K \boldsymbol{n}_{K,\sigma} \cdot \boldsymbol{n}_{K,\sigma}}{d_{K,\sigma}^2} \psi(s)^2 [\vartheta]_{K,\sigma}(\boldsymbol{y}) [\phi]_{K,\sigma}(\boldsymbol{y}) \, d\boldsymbol{x}$$

$$= I + II + III.$$

Since $\boldsymbol{x} = \boldsymbol{x}_K + s(\boldsymbol{y} - \boldsymbol{x}_K)$ for $s \in ]0,1[$, $\boldsymbol{y} \in \sigma$ and $\nabla_{\overline{K}} \vartheta \in \mathbb{P}_{\ell-1}(K)^d$ then

$$\nabla_{\overline{K}} \vartheta(\boldsymbol{x}) \cdot \boldsymbol{n}_{K,\sigma} = \nabla_{\overline{K}} \vartheta(\boldsymbol{y}) \cdot \boldsymbol{n}_{K,\sigma} + \sum_{j=1}^{\ell-1} p_j(\boldsymbol{y})(1-s)^j,$$

with $p_j(\boldsymbol{y})$ polynomials of $\ell - 1$ degree in the components of $\boldsymbol{y}$. It follows from Eq. (2.7) that

$$\int_0^1 \nabla_{\overline{K}} \vartheta(\boldsymbol{x}) \cdot \boldsymbol{n}_{K,\sigma} s^{d-1} \psi(s) \, ds = \nabla_{\overline{K}} \vartheta(\boldsymbol{y}) \cdot \boldsymbol{n}_{K,\sigma},$$

hence, using the change of variables $d\boldsymbol{x} = s^{d-1} d_{K,\sigma} ds d\boldsymbol{y}$ we get

$$II = \sum_{K \in \mathcal{M}} \sum_{\sigma \in \mathcal{F}_K} \int_\sigma A_K ([\vartheta]_{K,\sigma}(\boldsymbol{y}) \nabla_{\overline{K}} \phi(\boldsymbol{y}) + [\phi]_{K,\sigma}(\boldsymbol{y}) \nabla_{\overline{K}} \vartheta(\boldsymbol{y})) \cdot \boldsymbol{n}_{K,\sigma} \, d\boldsymbol{y}.$$

For $\sigma \in \mathcal{F}_i$ with $\sigma = \partial K \cap \partial T$ let $\boldsymbol{n}_\sigma = \boldsymbol{n}_{K,\sigma}$ and

$$[\![\Pi_\mathcal{D} \vartheta]\!]_\sigma = \Pi_{\overline{K}} \vartheta - \Pi_{\overline{T}} \vartheta, \qquad \{\!\!\{ A \nabla \Pi_\mathcal{D} \vartheta \}\!\!\}_{\omega,\sigma} = \omega_{K,\sigma} A|_K \nabla \Pi_{\overline{K}} \vartheta + \omega_{T,\sigma} A|_T \nabla \Pi_{\overline{T}} \vartheta.$$

If $\sigma \in \mathcal{F}_b$ with $\sigma = \partial K \cap \partial \Omega$ let $\boldsymbol{n}_\sigma = \boldsymbol{n}_{K,\sigma}$ and

$$[\![\Pi_\mathcal{D}\vartheta]\!]_\sigma = \Pi_{\overline{K}}\vartheta, \qquad \{\!\{A\nabla\Pi_\mathcal{D}\vartheta\}\!\}_{\omega,\sigma} = A|_K \nabla\Pi_{\overline{K}}\vartheta.$$

It holds $[\vartheta]_{K,\sigma} \cdot \boldsymbol{n}_{K,\sigma} = -\omega_{K,\sigma}[\![\Pi_\mathcal{D}\vartheta]\!]_\sigma \cdot \boldsymbol{n}_\sigma$ and similarly for $\phi$, hence

$$II = -\sum_{K\in\mathcal{M}} \sum_{\sigma\in\mathcal{F}_K} \int_\sigma \omega_{K,\sigma}([\![\Pi_\mathcal{D}\vartheta]\!]_\sigma A_K \nabla\Pi_{\overline{K}}\phi + [\![\Pi_\mathcal{D}\phi]\!]_\sigma A_K \nabla\Pi_{\overline{K}}\vartheta) \cdot \boldsymbol{n}_\sigma \, d\boldsymbol{y}$$

$$= -\sum_{\sigma\in\mathcal{F}} \int_\sigma ([\![\Pi_\mathcal{D}\vartheta]\!]_\sigma \{\!\{A\nabla\Pi_\mathcal{D}\phi\}\!\}_{\omega,\sigma} + [\![\Pi_\mathcal{D}\phi]\!]_\sigma \{\!\{A\nabla\Pi_\mathcal{D}\vartheta\}\!\}_{\omega,\sigma}) \cdot \boldsymbol{n}_\sigma \, d\boldsymbol{y}.$$

For $III$, using $C_\psi^2 = \int_\alpha^1 \psi(s)^2 s^{d-1} ds$ and by the usual change of variables, we obtain

$$III = C_\psi^2 \sum_{K\in\mathcal{M}} \sum_{\sigma\in\mathcal{F}_K} \frac{\delta_{K,\sigma}}{d_{K,\sigma}} \omega_{K,\sigma}^2 \int_\sigma [\![\Pi_\mathcal{D}\vartheta]\!]_\sigma [\![\Pi_\mathcal{D}\phi]\!]_\sigma \, d\boldsymbol{y}$$

$$= \sum_{\sigma\in\mathcal{F}} \eta_\sigma \frac{\gamma_\sigma}{h_\sigma} \int_\sigma [\![\Pi_\mathcal{D}\vartheta]\!]_\sigma [\![\Pi_\mathcal{D}\phi]\!]_\sigma \, d\boldsymbol{y},$$

where $h_\sigma$ is the diameter of $\sigma$ and $\gamma_\sigma$, $\eta_\sigma$ for $\sigma \in \mathcal{F}_i$ are defined by

$$\gamma_\sigma = \frac{2\delta_{K,\sigma}\delta_{T,\sigma}}{\delta_{K,\sigma} + \delta_{T,\sigma}},$$

$$\eta_\sigma = C_\psi^2 \left( \frac{\delta_{K,\sigma}}{d_{K,\sigma}} \omega_{K,\sigma}^2 + \frac{\delta_{T,\sigma}}{d_{T,\sigma}} \omega_{T,\sigma}^2 \right) \frac{h_\sigma}{\gamma_\sigma} = C_\psi^2 h_\sigma \left( \frac{\omega_{K,\sigma}}{d_{K,\sigma}} + \frac{\omega_{T,\sigma}}{d_{T,\sigma}} \right)$$

and for $\sigma \in \mathcal{F}_b$ by

$$\gamma_\sigma = \delta_{K,\sigma}, \qquad\qquad \eta_\sigma = C_\psi^2 \frac{h_\sigma}{d_{K,\sigma}}.$$

Summing $I, II, III$ we get the equivalence of $\int_\Omega A\nabla_\mathcal{D}\vartheta \cdot \nabla_\mathcal{D}\phi \, d\boldsymbol{x}$ and [7, equation 4.64] with the parameter $\eta_\sigma$ chosen as above. Under the additional hypothesis that the mesh sequence satisfies

$$\min\{\frac{h_\sigma}{d_{K,\sigma}} : K \in \mathcal{M}, \sigma \in \mathcal{F}_K\} \geq C_\mathcal{F} > 0,$$

we have that $\eta_\sigma \geq C_\psi^2 C_\mathcal{F}$. Since $C_\psi^2 \geq d/(1-\alpha^d)$, letting $\alpha \to 1$ we can have $\eta_\sigma$ as large as desired.